\documentclass{amsart}
\newtheorem{thm}{Theorem}[section]

\newtheorem{defn}[thm]{Definition}
\newtheorem{lemma}[thm]{Lemma}
\newtheorem{conj}[thm]{Conjecture}
\newtheorem{cor}[thm]{Corollary}
\newtheorem{remark}[thm]{Remark}

\usepackage{amsmath}
\usepackage{amsxtra}
\usepackage{amscd}
\usepackage{amsthm}
\usepackage{amsfonts}
\usepackage{amssymb}
\usepackage{eucal}
\newcommand{\cQ}{\mathcal Q}
\newcommand{\cT}{\mathcal T}
\newcommand{\cF}{\mathcal F}
\newcommand{\cP}{{\rm PS}}

\newcommand{\g}{{\mathfrak{g}}}
\newcommand{\h}{{\mathfrak{h}}}

\newcommand{\n}{{\mathfrak{n}}}
\renewcommand{\sl}{{\mathfrak{sl}}}
\newcommand{\KR}{{{\rm KR}}}
\newcommand{\ch}{{\rm ch}}
\newcommand{\C}{{\mathbb C}}
\newcommand{\Z}{{\mathbb Z}}
\newcommand{\N}{{\mathbb N}}

\newcommand{\bm}{{\mathbf m}}
\newcommand{\bn}{{\mathbf n}}
\newcommand{\bu}{{\mathbf u}}
\newcommand{\al}{{\alpha}}

\numberwithin{equation}{section}

\begin{document}
\title[Combinatorial Kirillov-Reshetikhin conjecture]{Proof of the combinatorial Kirillov-Reshetikhin conjecture}

\author{P. Di Francesco and R. Kedem}
\address{PDF: Service de Physique Th\'eorique de Saclay, CEA/DSM/SPhT,
  URA 2306 du CNRS, F-91191 Gif sur Yvette Cedex.}
\address{RK: Department of Mathematics, University of Illinois, 1409 W
  Green Street, Urbana IL 61801 USA.}

\begin{abstract}
In this paper we give a direct proof of the equality of certain generating
function associated with tensor product multiplicities of
Kirillov-Reshetikhin modules of the untwisted Yangian for each simple Lie algebra
$\g$. Together with the theorems of Nakajima and Hernandez, this
gives the proof of the combinatorial version of the Kirillov-Reshetikhin
conjecture, which gives tensor product multiplicities in terms of
restricted fermionic summations.
\end{abstract}

\maketitle

\section{Introduction}
This paper aims to resolve three related conjectures. As explained below, the
proof of one implies the proof of the other two. Here, we introduce a general method which allows us to prove
the "$M=N$ conjecture" of \cite{HKOTY}. A brief sketch of the
necessary conjectures and theorems follows.

The Kirillov-Reshetikhin conjecture about completeness of Bethe states
in the generalized inhomogeneous Heisenberg spin-chain is a
combinatorial formula (``the $M$-sum'' \eqref{M}) for the number of solutions of the Bethe equations. The formula has a fermionic form, in that it is a sum over
products of non-negative binomial coefficients. We call this
 the
combinatorial Kirillov-Reshetikhin conjecture, Conjecture
\ref{KRconjecture}.

This is closely related to another conjecture, now a theorem (Theorem \ref{NakHer}), also
sometimes referred to in the literature as the Kirillov-Reshetikhin
conjecture, about the characters of special finite-dimensional quantum
affine algebra modules.  This second version was recently proved by
Nakajima \cite{Nakajima} for simply-laced Lie algebras, and more
generally by Hernandez \cite{Hernandez}. It is the statement that
characters of Kirillov-Reshetikhin modules are solutions of a
recursion relation called the $Q$-system \cite{KR,HKOTY}. This has
been proven to hold for all simple Lie algebras, and also in more
general settings \cite{He07}.

Hatayama et. al. \cite{HKOTY} proved that this theorem implies a
certain explicit alternating sum formula in terms of binomial coefficients
(``the $N$-sum'' \eqref{N}) for the tensor product multiplicities. This formula
is closely related to, but manifestly different from, the $M$-sum in
the combinatorial Kirillov-Reshetikhin conjecture, which involves a
restricted, non-alternating sum.

The second conjecture, then, is the one advanced in  \cite{HKOTY}, that the two sums are
equal. We call this the ``$M=N$ conjecture'' (the precise statement is
Conjecture \ref{HKOTY}). There is an indirect argument which shows
that this is true in special cases, because the $M$-sum
formula was proven by combinatorial means in certain special cases by
\cite{KKR,KR,KS,KSS,OkaSchShi,Schilling}. In this paper we prove this
conjecture directly, thereby proving the combinatorial
Kirillov-Reshetikhin conjecture.

The third conjecture is the Feigin-Loktev conjecture. 
It was shown in \cite{AK} that the combinatorial KR-conjecture implies
the Feigin-Loktev conjecture \cite{FL} for the fusion product of
arbitrary KR-modules of $\g[t]$, as defined by \cite{Chari,ChariMoura}. This
conjecture states that the graded tensor product (the unrestricted
Feigin-Loktev fusion product) is independent of the localization
parameters, and that its graded dimension is given by the generalized
Kostka polynomials \cite{SS} or fermionic sums of \cite{HKOTY}.

In this paper, we prove the $M=N$ conjecture in the untwisted case for any simple Lie algebra, by considering suitable
generating functions, a standard technique in enumerative
combinatorics.  These functions are constructed so as to enjoy
particularly nice factorization properties, in terms of solutions of
the so-called $Q$-systems or certain deformations thereof.  Analogous
generating functions, involving fewer parameters, were also used in
\cite{HKOTY} in this context, but for a different purpose.

The proof of Conjecture \ref{HKOTY} completes the proof of the combinatorial
 Kirillov-Reshetikhin conjecture for representations of $Y(\g)$ for
 all simple Lie algebras $\g$ (we do not consider twisted cases in
 this paper). 

In addition, in the cases where Chari's
 KR-modules for $\g[t]$ are known to have the same dimension as their
 Yangian version (the classical algebras and some exceptional cases), 
this also completes the proof of the Feigin-Loktev
 conjecture for the unrestricted fusion products of arbitrary
 Kirillov-Reshetikhin modules for any simple Lie algebra.

In the course of our proof, we define a new family of functions,
generalizing the characters $Q$ of the KR-modules, which appear to
have very useful properties. First, we define a deformation of the
$Q$-system, in terms of functions in an increasing number of
variables.  We show that these functions can be defined alternatively
in terms of a substitution recursion. In terms of the deformed
$Q$-functions, there is a complete factorization of generating
functions for fermionic sums of the KR-type.

The paper is organized as follows. In Section 2, we recall the
definitions, conjectures and theorems which we use in this paper. 
In Sections 3, 4, and 5, we give the proof of the $M=N$ conjecture for
for $\g=\sl_2$,
$\g$ simply-laced and $\g$ non simply-laced respectively. In each
case, we define a deformed $Q$-system, which we refer to as the
$\cQ$-system. We then define generating functions of fermionic sums
which have a factorized form in terms of the $\cQ$-functions. This
allows us to prove an equality of restricted generating functions, and
the constant term of this identity is the conjectured $M=N$
identity of \cite{HKOTY}.

\section{$Q$-systems, the Kirillov-Reshetikhin conjecture and the Feigin-Loktev conjecture}
\subsection{Definitions}
Let $\g$ be a simple Lie algebra with simple roots $\alpha_i$ with
$i\in I_r = \{1,...,r\}$ and Cartan matrix $C$ with entries $C_{i,j} =
\frac{2 (\alpha_i,\alpha_j)}{(
\alpha_i,\alpha_i)}$. The algebra $\g$ has Cartan decomposition
$\g=\n_- \oplus \h \oplus \n_+$, and we denote the generator of $\n_-$
corresponding to the simple root $\al_i$ by $f_i$ etc..

The irreducible integrable highest weight modules of $\g$ are denoted
by $V(\lambda)$, where $\lambda\in P^+$ and $P^+$ is the set of
dominant integral weights. We denote the fundamental weights of $\g$ by
$\omega_i$ ($i\in I_r$).

The algebra $\g[t] = \g\otimes \C[t]$ is the Lie algebra of
polynomials in $t$ with coefficients in $\g$. The generators of
$\g[t]$ are denoted by $x[n]:=x\otimes t^n$ where $x\in \g$ and $n\in
\Z_+$. The relations in the algebra are
$$
[x\otimes f(t),y\otimes h(t)]_{\g[t]} = [x,y]_\g f(t) h(t), \qquad x,y\in
\g,\ f(t),h(t)\in \C[t],
$$ 
where $[x,y]_\g$ is the usual Lie bracket in $\g$.
We regard $\g$ as the subalgebra of constant currents in $\g[t]$. Thus,
any $\g[t]$-module is also a $\g$-module by restriction.

\subsubsection{Localization}
Any $\g$-module $V$ can be extended to a $\g[t]$-module by the
evaluation homomorphism. That is, given a complex number $\zeta$, the
module $V(\zeta)$ is the $\g[t]$-module defined by
$$
x[n] v = \zeta^n x v,\quad v\in V(\zeta).
$$ If $V$ is finite-dimensional, so is $V(\zeta)$, with the same
dimension. If $V$ is irreducible as a $\g$-module, so is $V(\zeta)$.

More generally, given a $\g[t]$-module $V$, the $\g[t]$-module
localized at $\zeta$, $V(\zeta)$, is the module on which $\g[t]$ acts
by expansion in the local parameter $t_\zeta := t-\zeta$. If $v\in
V(\zeta)$, then
$$x[n] v = x\otimes (t_\zeta + \zeta)^n v = \sum_j {n \choose j}
\zeta^j x[n-j]_\zeta v,$$ where $x[n]_\zeta := x\otimes t_\zeta$ and
$x[n]_\zeta$ acts on $v\in V(\zeta)$ in the same way that $x[n]$ acts
on $v\in V$.

An evaluation module $V(\zeta)$ is a special case of a localized
module, on which the positive modes $x[n]_\zeta$ with $n>0$ and
$x\in\g$ act
trivially.

\subsubsection{Grading}\label{gradingsection}
Let $V$ be any cyclic $\g[t]$-module. That is, there exists a vector
$v\in V$
such that 
$$
V = U(\g[t]) v.
$$
Any such module can be endowed with a $\g$-equivariant grading (which
depends on the choice of $v$ in case it is non-unique) as
follows. The algebra $U=U(\g[t])$ is graded by degree in $t$. That is,
the graded component $U^{(j)}$ is the span of monomials of the form
$$
x_1[n_1]\cdots x_m[n_m]:\quad \sum_{i=1}^m n_i = j,
$$ where $x_i\in\g$ and $m\in \Z_+$. 

The module $V$ does not necessarily inherit this grading since the
action of $\g[t]$ is not assumed to respect this grading. This is
true, in general, for the localized modules above if $\zeta\neq
0$. However, $V$ does inherit a filtration, which depends on the
choice of $v$. Let $U^{(\leq i)}$ be the vector space generated by
monomials with degree less than or equal to $i$ in $U$. Define $\cF{(i)}
= U^{(\leq i)} v$.  We have $\cF{(i)}\subset \cF{(i+1)}$, where $\cF{(0)}$
is the $\g$-module generated by $v$.

In the case that $V$ is finite-dimensional, this gives a finite
filtration of $V$,
$$
\cF{(0)}\subset \cdots \subset \cF{(N)} = V.
$$
The associated graded space (here $\cF{(-1)}=\emptyset$),
$$
{\rm Gr}\ V = \underset{i\geq 0}{\oplus} \cF{(i)}/\cF{(i-1)}
$$
has graded components $V[i] = \cF{(i)}/\cF{(i-1)}$ which are
$\g$-modules, since the filtration is $\g$-equivariant.

\subsection{Kirillov-Reshetikhin modules}
The term Kirillov-Reshetikhin module properly refers to certain
finite-dimensional Yangian modules \cite{KR} or quantum affine algebra
modules.  Chari's Kirillov-Reshetikhin modules \cite{Chari,ChariMoura} are
$\g[t]$-modules which are classical limits of the quantum group
modules. Whereas $Y(\g)$ and $U_q(\widehat{\g})$-modules are defined
in terms of their Drinfeld polynomial, Chari's KR-modules for $\g[t]$
are defined in terms of generators and relations.  We refer to Chari's
modules as KR-modules in this paper. As $\g$-modules, they are known
to have the same structure as the Yangian KR-modules in the case of
the classical Lie algebras, and in certain exceptional cases. 

These modules arise naturally when one considers the explicit
description of the dual space of functions \cite{AK} to Feigin-Loktev
fusion product \cite{FL}.

KR-modules are parametrized by a complex number $\zeta\in \C^*$ (the
localization parameter) and a highest weight of the special form $m
\omega_\al\ (\al\in I_r)$, with $m\in \Z_+$ and $\omega_\al$ a simple
weight. We denote such a module by $\KR_{\al,m}(\zeta)$. The
definition given here is the one used in \cite{AK}.

\begin{defn}
Let $\zeta\in C^*$ and let $m\in \Z_+$, $\al\in I_r$. The {\rm KR}-module
$\KR_{\al,m}(\zeta)$ is the module generated by the action of $U(\g[t])$
on a cyclic vector $v\in \KR_{\al,m}(\zeta)$, subject to the relations
(recall that $x[n]_\zeta = x\otimes (t-\zeta)^n$):
\begin{eqnarray*}
x[n]_\zeta v &=& 0 \quad \hbox{if $x\in \n_+$ and $n\geq 0$};\\
h_\beta[n]_\zeta v &=& \delta_{n,0} \delta_{\al,\beta} m v;\\ 
f_\beta[n]_\zeta v &=& 0 \quad \hbox{if $n\geq \delta_{\al,\beta}$};\\
f_\al[0]_\zeta^{m+1} v &=& 0.
\end{eqnarray*}
The associated graded space of this module is a graded $\g[t]$-module
$\overline{\KR}_{\al,m}$. Its graded components are $\g$-modules.
\end{defn}

For example, in the case of $\g=A_r$, $\KR_{\al,m}(\zeta) =
V_{m\omega_\al}(\zeta)$, the evaluation module of $\sl_{r+1}[t]$ corresponding
to the irreducible $\g$-module $V(m \omega_\al)$ at the point $\zeta$.
For other Lie algebras, $\KR_{\al,m}(\zeta)$ may not be irreducible as
a $\g$-module. However, 
the decomposition of $\KR_{\al,m}(\zeta)$ into irreducible $\g$-modules
always has a unitriangular form (in the partial ordering of weights). That is,
$$
\KR_{\al,m}(\zeta) \underset{\g-{\rm mod}}{\simeq}
V(m\omega_\al)\oplus \left(\underset{\mu<m\omega_\al}{\oplus}
V(\mu)^{\oplus m_{\mu}}\right)
$$
Thus, $m\omega_\al$ is the highest $\g$-weight of $\KR_{\al,m}(\zeta)$.

The decomposition of tensor products of KR-modules
into irreducible $\g$-modules, is the subject of the
Kirillov-Reshetikhin conjecture.

\subsection{The Kirillov-Reshetikhin conjecture}

The KR-conjecture was originally a conjecture about the completeness
of Bethe ansatz states for the generalized, inhomogeneous Heisenberg
spin chain. This is a spin-chain model with inhomegeneity parameters
$\zeta_i$ at each lattice site $i$, and with a representation
$V_i(\zeta_i)$ of the Yangian $Y(\g)$ (or, equivalently, of
$U_q(\widehat{\g})$ if the generalized XXZ-model is considered) at
each lattice site. 

The modules $V_i(\zeta_i)$ are each assumed to be of
Kirillov-Reshetikhin type. If the Bethe ansatz gives a complete set of
solutions, then the Bethe states should be in one-to-one
correspondence with $\g$-highest weight vectors in the Hilbert space
of the Hamiltonian or transfer matrix. The Hilbert space is simply the tensor
product of the modules $V_i(\zeta_i)$.

The Bethe states are parametrized by solutions of certain coupled
algebraic equations, known as the Bethe equations. It is hypothesized
that the solutions of the Bethe ansatz equations have a certain form
of their complex parts, the so-called ``string hypothesis,'' and, more
importantly in this context, are parametrized by the Bethe integers.

\begin{remark}
It is known that the string hypothesis is not, in fact, correct in
general. However, solutions to the Bethe equations can still be shown
to be parametrized by the Bethe integers in certain cases. The
correctness of the string hypothesis is not relevant for the current
paper. It served only as the inspiration for the original
Kirillov-Reshetikhin conjecture.
\end{remark}

The Kirillov-Reshetikhin conjecture is that the Bethe integers
parametrize solutions of the Bethe equations. It can be formulated in
completely combinatorial terms as follows.

Let $\bn = \{n_{\al,i} |\ \al\in I_r, i\in \N\}$ be a collection of
non-negative integers whose sum is finite.
These parametrize a set of $N=\sum_{i,\al}
n_{\al,i}$ KR-modules, with $n_{\al,i}$ modules with highest
$\g$-weight $i\omega_\al$, and hence they parametrize the Hilbert space.

For each $\lambda$ a dominant integral weight
$$\lambda= \sum_{\al\in I_r} l_\al \omega_\al\in P^+ $$
choose
a set of non-negative integers
$\bm=\{m_{\al,i}\}$ with $\al\in I_r$ and $i\in\N$, such that the total spin
\begin{equation}\label{spin}
q_\al = l_\al + \sum_{i,\beta} i C_{\al,\beta} m_{\beta,i} - \sum_i
i n_{\al,i}, \quad (\al\in I_r)
\end{equation}
is zero, $q_\al=0$.

Define the ``vacancy numbers'' which depend on the sets of integers
$\{l_\al\}, \bm, \bn$ and on the Cartan matrix:
\begin{equation}\label{pvacancy}
p_{\al,i} = \sum_{j\geq 1} n_{\al,j} \min(i,j) - \sum_{\beta\in I_r}{\rm
sgn}(C_{\al,\beta}) \sum_{j\geq 1} \min(|C_{\al,\beta}| j, |C_{\beta,\al}|i)
m_{\beta,j}, \quad(\al\in I_r, i\in \N).
\end{equation}

The Bethe integers are any set of $m_{\al,i}$ distinct integers chosen
from the interval $[0,p_{\al,i}]$ for each $\al$ and $i$. Therefore
$p_{\al,i}<0$ does not correspond to any Bethe states. The number of
distinct sets of Bethe integers is The fermionic multiplicity formula
called the $M$-sum:
\begin{equation}\label{M}
M_{\lambda;\bn} = \sum_{\underset{q_\al=0,p_{\al,i}\geq
    0}{m_{\al,i}\geq 0}} \prod_\al  
{m_{\al,i}+p_{\al,i}\choose m_{\al,i}}.
\end{equation}
Here, the summation is over all non-negative integers $\{m_{\al,i}\}$.
For fixed $\bn$ and $\lambda$ there is some integer $p$ such
that all $m_{\al,i}$ with $i>p$ are constrained to be zero (due to the
constraint $q_\al=0$), so that there is only a finite number of
summation variables.

\begin{remark}
One can attach an ``energy'' to each Bethe integer which is proportional to
the integer itself. In this way, one obtains a graded multiplicity
formula $M_{\lambda;\bn}(q)$ which is a polynomial keeping track of
the energy grading parameter. This grading appears also in the fusion
product, described below. Although it is of interest in discussing the
fusion product, for the proof of the identities in this paper, it is
not necessary to keep track of this grading.
\end{remark}

Note that the binomial coefficients are defined for both positive and
negative values of $p_{\al,i}$:
$$
{m+p\choose m} = \frac{(p+m)(p+m-1)\cdots (p+1)}{m!}.
$$ If $p<0$ then if $m<-p$, the sign of the binomial coefficient is
$(-1)^m$. In general, the summation over the variables $m_{\al,i}$
might include both negative and positive terms. In the $M$-sum, terms with
$p_{\al,i}<0$ are excluded.

\begin{conj}[The combinatorial Kirillov-Reshetikhin conjecture
    \cite{KR,HKOTY}]\label{KRconjecture}
\begin{equation}
{\rm dim}\ {\rm Hom}_\g\ \left( \underset{\al,i}{\otimes} \KR_{\al,i}^{\otimes n_{\al,i}} ,\
V(\lambda)\right) = M_{\lambda,\bn}.
\end{equation}
\end{conj}

This conjecture has been proven for the following special cases of
$\g$ and $\bn$:
\begin{itemize}
\item For $\g=A_r$, the conjecture was proven by \cite{KKR,KR} and
  \cite{KSS} for arbitrary $\bn$.
\item For $\g=D_r$, the conjecture was proven in \cite{Schilling} for
  $\bn$ such that $n_{\al,1}\geq 0$ and $n_{\al,j}=0$ for all
  $j>1$.
\item For $\g$ any non-exceptional simple Lie algebra and $\bn$ such
  that $n_{1,i}\geq 0$ and $n_{\al,j}=0$ for all $\al>1$ \cite{OkaSchShi,SS}.
\end{itemize}
The proof in each of these cases involves a bijection between 
combinatorial objects known as ``rigged configurations'' and crystal
paths.

\subsection{$Q$-systems}\label{q-system}

Kirillov and Reshetikhin \cite{KR} also introduced another, closely
related conjecture, recently proven for all simple Lie algebras $\g$
\cite{Nakajima,Hernandez} and some generalizations \cite{He07} (more
precisely, the conjecture is concerned with finite-dimensional modules
of $U_q(\widehat{\g})$ or $Y(\g)$).

\begin{thm}[\cite{KR,Nakajima,Hernandez}]\label{NakHer}
The characters $Q_{\al,i}$ of the Kirillov-Reshetikhin modules of
$U_q(\widehat{\g})$ for any simple Lie algebra $\g$ satisfy the
so-called $Q$-system (Equation \eqref{qsys} below). In addition, they
satisfy \cite{Hernandez} the asymptotic conditions of \cite{HKOTY} (condition C of Theorem 7.1 \cite{HKOTY}), so
that their decomposition into irreducible $U_q(\g)$-modules is given
by Equation \eqref{N}.
\end{thm}

The $Q$-system is a quadratic recursion relation for the the family of
functions $\{Q_{\al,j}: \al\in I_r, j\in \N\}$. Each element
$Q_{\al,j}$ has the interpretation of the character of the
Kirillov-Reshetikhin module corresponding to a highest $\g$-weight
$j\omega_\al$, where $\omega_\al$ is one of the fundamental weights of
$\g$.

In general the recursion relation is \cite{KR,HKOTY}
\begin{equation}\label{qsys}
Q_{\al,j+1} = \frac{\displaystyle Q_{\al,j}^2-\prod_{\beta\sim
\al}\prod_{k=0}^{|C_{\al,\beta}|-1}
Q_{\beta,\lfloor(|C_{\beta,\al}|j+k)/|C_{\al,\beta}|\rfloor}}{Q_{\al,j-1}},
\quad
(j>0),
\end{equation}
with initial conditions $Q_{\al,0}=1$ and $Q_{\al,1}=t_\al$, a formal
variable.
Here, $\beta\sim\al$ means that the nodes $\al$ and $\beta$ are
connected in the Dynkin diagram of $\g$. The notation $\lfloor a
\rfloor$ denotes the integer part of $a$.

If $\g$ is a simply-laced Lie algebra, then the system has the form
\begin{equation}
Q_{\al,j+1} = \frac{Q_{\al,j}^2 - \prod_{\beta\sim \al}
Q_{\beta,j} }{Q_{\al,j-1}}, \quad (j>0).
\end{equation}

In the non-simply laced case, the relations have the form
\begin{equation}
Q_{\al,j+1} = \frac{Q_{\al,j}^2 - \prod_{\beta\sim \al}
T^{(\al,\beta)}_j }{Q_{\al,j-1}},
\end{equation}
where $T_j^{(\al,\beta)}=Q_{\beta,j}^{|C_{\al,\beta}|}$ except in the
following cases:
\subsection*{$B_r$}
\begin{eqnarray*}
T_j^{(r-1,r)} &=& Q_{r,2j}\\
T_j^{(r,r-1)} &=& Q_{r-1,\lfloor j/2\rfloor}Q_{r-1,\lfloor (j+1)/2\rfloor}.
\end{eqnarray*}
\subsection*{$C_r$}
\begin{eqnarray*}
T_j^{(r-1,r)} &=& Q_{r,\lfloor j/2\rfloor}Q_{r,\lfloor (j+1)/2\rfloor}\\
T_j^{(r,r-1)} &=& Q_{r-1,2j}.
\end{eqnarray*}
\subsection*{$F_4$}
\begin{eqnarray*}
T_j^{(3,2)} &=& Q_{2,\lfloor j/2\rfloor}Q_{2,\lfloor (j+1)/2\rfloor}\\
T_j^{(2,3)} &=& Q_{3,2j}.
\end{eqnarray*}
\subsection*{$G_2$}
\begin{eqnarray*}
T_j^{(2,1)} &=& Q_{1,\lfloor j/3\rfloor}Q_{1,\lfloor
(j+1)/3\rfloor}Q_{1,\lfloor (j+2)/3\rfloor}\\ 
T_j^{(1,2)} &=&
Q_{2,3j}.
\end{eqnarray*}

We note an important corollary of this fact, which we call the
polynomiality property of KR-characters:
\begin{thm}\label{polynomiality}
Given the data $\{Q_{\al,0}=1\}_{\al\in I_r}$, the solutions of the
$Q$-system are polynomials in the variables $\{Q_{\al,1}\}_{\al\in I_r}$.
\end{thm}
\begin{proof}
This is simply the statement that the Groethendieck group of
KR-modules is generated by the trivial representation and the
fundamental KR-module with highest weight $\omega_\al$.

It is known that there is a unitriangular
decomposition of the KR-characters into irreducible
$U_q(\g)$-characters, with the highest weight module $V(i \omega_\al)$
appearing with multiplicity one in $KR_{\al,i}$. The other modules in
the decomposition have highest weights which are strictly lower in the
partial ordering of weights. The statement of polynomiality follows
from this fact.
\end{proof}

It was proved in \cite{HKOTY} that if the characters $\{Q_{\al,i}\}$
satisfy the $Q$-system plus a certain asymptotic condition, then the
characters of their tensor products have an explicit ``fermionic''
expression
\begin{thm}[Theorem 8.1 \cite{HKOTY}]\label{HKOTYthm}
Define $Q_{\al,i}$ to be the $U_q(\g)$-character of the KR-module
corresponding to highest weight $i \omega_\al$. Then
\begin{equation}
\prod_{\al,i} Q_{\al,i}^{n_{\al,i}} = \sum_{\lambda} N_{\lambda;\bn}
\ch V(\lambda),
\end{equation}
where $V(\lambda)$ is the irreducible $U_q(\g)$-module with highest weight
$\lambda$. 
\end{thm}

Here, the $N$-sum is
\begin{equation}\label{N}
N_{\lambda,\bn} = \sum_{m_{\al,i}\geq 0\atop q_\al=0 } \prod_\al  
{m_{\al,i}+p_{\al,i}\choose m_{\al,i}},
\end{equation}
where $q_\al$ and $p_{\al,i}$ are defined by \eqref{spin} and
\eqref{pvacancy} as for the $M$-sum. The only difference between this
conjecture and the combinatorial KR-conjecture \ref{KRconjecture} is
that the summation is not restricted to non-negative values of the
vacancy numbers. That is, the $N$-sum has more terms, some of which
are negative.

In fact, \cite{HKOTY} conjectured that the two sums are equal. We will
describe something which we call the HKOTY-conjecture, which is
slightly stronger than this. (Their conjecture extends to the graded
dimensions, which we will introduce below for fusion products. However,
we need only prove the following version.)  

Define $N_{\lambda,\bn}^{(k)}$ and $M_{\lambda,\bn}^{(k)}$ to be the
sums in equations \eqref{N} and \eqref{M}, respectively, with the
summations further restricted so that $m_{\al, i}=n_{\al,i}=0$ if $i>t_\al
k$. (Here, $t_\al$ is 1 for the long roots, $2$ for the short
roots of $B_r, C_r, F_4$ and 3 for the short root of $G_2$.)
\begin{conj}[The HKOTY-conjecture]\label{HKOTYconj}
For any simple Lie algebra,
\begin{equation}\label{mainidentity}
M_{\lambda,\bn}^{(k)} = N_{\lambda,\bn}^{(k)}.
\end{equation}
\end{conj}

We have $M_{\lambda,\bn} = \lim_{k\to\infty} M_{\lambda,\bn}^{(k)}$
and $N_{\lambda,\bn} = \lim_{k\to\infty}
N_{\lambda,\bn}^{(k)}$. Therefore, if the conjecture \ref{HKOTYconj}
is true, then combined with Theorem \ref{HKOTYthm} and the result of
\cite{Hernandez}, it implies the completeness conjecture
\ref{KRconjecture} of Kirillov and Reshetikhin.

The purpose of the current article is to prove this conjecture
directly, for all simple Lie algebras and for all $\bn$.

\subsection{The Feigin-Loktev conjecture and Kirillov-Reshetikhin conjecture}

In their paper \cite{FL} the authors introduced a graded tensor
product on finite-dimensional, graded, cyclic $\g[t]$-modules for $\g$
a simple Lie algebra, which they call the fusion product. This is
related to the fusion product in Wess-Zumino-Witten conformal field
theory when the level is restricted. In the current paper, we take the
level to be sufficiently large so that it does not enter the
calculations. This is called the unrestricted fusion product.

Let us summarize the results of \cite{AK} concerning the
Feigin-Loktev conjecture for unrestricted fusion products of
Kirillov-Reshetikhin modules for any simple Lie algebra $\g$.
The description below of the Feigin-Loktev fusion product \cite{FL} is 
the one given in \cite{Kedem,AK}. We refer the reader to those
articles for further details.

Let $\{\zeta_1,...,\zeta_N\}$ to be distinct complex numbers,
and choose $\{V_1(\zeta_1),...,V_N(\zeta_N)\}$ to be
finite-dimensional $\g[t]$-modules, cyclic with cyclic vectors $v_i$,
localized at the points $\zeta_i$. Then as $\g$-modules, we have
\cite{FL},
\begin{equation}\label{filtered}
V_1(\zeta_1)\otimes\cdots\otimes V_N(\zeta_N) \simeq U(\g[t])
v_1\otimes\cdots\otimes v_N.
\end{equation}

This tensor product is also a cyclic $\g[t]$-module with cyclic vector
$v_1\otimes \cdots \otimes v_N$. It inherits a filtration from
$U(\g[t])$ as in Section \ref{gradingsection}.

\begin{defn}
The unrestricted Feigin-Loktev fusion product, or graded tensor
product, is the associated graded space of the filtered space
\eqref{filtered}. It is denoted by $V_1\star \cdots
\star V_N (\zeta_1,...,\zeta_N)$.
\end{defn}

Note that the filtration, and hence the grading, is $\g$-equivariant,
and hence the graded components are $\g$-modules.  Let
$M_{\lambda;\{V_i\}}[n]$ denote the multiplicity of the irreducible
$\g$-module $V(\lambda)$ in the $n$th-graded component of the fusion
product $V_1\star\cdots\star V_N (\zeta_1,...,\zeta_N)$.

\begin{defn}
The graded multiplicity ($q$-multiplicty) of $V(\lambda)$ in
 the Feigin-Loktev fusion
product $V_1\star \cdots \star V_N (\zeta_1,...,\zeta_N)$ is the
 polynomial in $q$ defined as
$$
M_{\lambda;\{V_i\}} (q) = \sum_{n\geq 0} M_{\lambda;\{V_i\}} [n] q^n.
$$
\end{defn}

\begin{conj}[Feigin-Loktev \cite{FL}]\label{FLconj}
In the cases where $V_i$ are sufficiently well-behaved,
$M_{\lambda,\{V_i\}} (q)$ is independent of the localization
parameters $\zeta_i$. 
\end{conj}
At this time, it is not known what ``sufficiently well-behaved''
means in general, and this remains an open problem. In this paper we
consider KR-modules, which we prove satisfy the necessary criteria.

This in particular implies that the dimension of the fusion product
is equal to the dimension of the tensor product of the $\g$-modules
$V_i$, which is the $\g[t]$-module $V_i(\zeta_i)$ regarded as a $\g$-module.

We have a Lemma, which follows from the fact that the fusion product
 is a quotient of the filtered tensor product \eqref{filtered} and a
 standard deformation argument (Lemma 20 of \cite{FKLMM}),
 \cite{FJKLM}:
\begin{lemma}\cite{FJKLM}\label{dimensions}
\begin{equation}
M_{\lambda;\{V_i\}}(1) \geq {\rm Dim}\left( {\rm Hom}_\g
\left(\underset{i}\otimes V_i,\ V(\lambda)\right)\right).
\end{equation}
\end{lemma}

One way of proving Conjecture \ref{FLconj} is to compute
$M_{\lambda,\{V_i\}}(q)$ explicitly, and to show that the polynomial is
independent of $\zeta_i$. This is, of course, a stronger result. We
have the following conjecture, inspired by
\cite{FL} and generalized and partially proven in
\cite{Kedem,FJKLM2,AKS,AK}. 

\begin{conj}\label{strongFLconj}
 In the case where $V_i$ are all of
  Kirillov-Reshetikhin-Chari type, the graded multiplicities
  $M_{\lambda,\{V_i\}}(q)$ are equal to the generalized Kostka
  polynomials or the fermionic sums $M_{\lambda,\bn}(q)$ of \cite{KR,HKOTY}.
\end{conj}
This conjecture implies \ref{FLconj} for these cases, because the
polynomials are indepenent of the localization parameters.

Various special cases of \ref{FLconj} have been proven. The case of
$\sl_2$ was proven in \cite{FJKLM2} by proving \ref{strongFLconj} (in
this case, the multiplicities are the usual co-charge Kostka
polynomials).  Conjecture \ref{FLconj} was proven for $\sl_n$
symmetric-power representations in \cite{Kedem} by using a result of
\cite{GP}. In
\cite{AKS}, we proved \ref{strongFLconj} for $\sl_n$ KR-modules by
using to a result of \cite{KSS} for the fermionic form of
generalized Kostka polynomials, which are the $q$-multiplicities in
the case of tensor products of KR-modules of $\sl_n$. Most generally, in
\cite{AK}, the following theorem was proven:

\begin{thm}[\cite{AK}]\label{ardonnekedem}
\begin{equation}
M_{\lambda;\{V_i\}} (q) \leq M_{\lambda;\bn}(q)
\end{equation}
where by the inequality we mean the inequality for the coefficients in
each power of $q$.
Here $n_{\al,i}$ is the number of KR-modules $\KR_{\al,i}$ in the
fusion product.
\end{thm}
Each of the coefficients on both sides is manifestly non-negative, so
it is sufficient to prove the equality for $q=1$, i.e. the equality of
total dimensions. Therefore, given Lemma \ref{dimensions}, in the
cases where Conjecture \ref{KRconjecture} has been proven, the set of
inequalities implies the equality of Hilbert polynomials, and hence
provides a proof of \ref{strongFLconj} and \ref{FLconj}.
Thus, in \cite{AKS,AK}, Conjecture \ref{strongFLconj} 
(hence \ref{FLconj}) was proven for these cases.

The proofs of \cite{KR,KSS,Schilling,OkaSchShi,SS} depend on a certain
 bijection from a combinatorial object called rigged configurations,
 which is what the $M$-summation counts, and crystal paths.  In the
 present paper, we bypass the question of the existence of crystal
 bases by proving the HKOTY conjecture directly. 

Together with Theorem 
 \ref{ardonnekedem}, this provides a proof of Conjectures
 \ref{strongFLconj}, \ref{FLconj} for all simple Lie algebras and the
 tensor products of any arbitrary set of KR-modules, modulo the
 identification of the dimension of Chari's KR-modules \cite{Chari}
 and the usual KR-modules. 

That is, Chari's modules (and the fusion
 product) are $\g[t]$-modules and not Yangian modules. It is known
 that the dimension of Chari's modules are equal to the dimensions of
 the KR modules for Yangians or, equivalently, quantum affine
 algebras, in the case of classical algebras and in some of the exceptional
 cases.

\section{Recursion relations and quadratic relations for $sl_2$}

In this section, we illustrate the method of the proof of the HKOTY
conjecture for the simplest case of $\sl_2$. The higher rank cases are
a straightforward generalization of this case, but the notation and formulas become much more
cumbersome. Hence it is instructive to examine this case first.

\subsection{The $\cQ$-system} 
We define a generalization of the $Q$-system for functions which we
call $\cQ$.  Let $\bu:=(u, u_i\ (i\in \N))$ be formal variables.  We
define a family of functions $\{\cQ_k(u;u_1,...,u_{k-1})\}_{k\in\Z_+}$
recursively as follows. It is convenient to use the notation $\cQ_k(\bu) :=
\cQ_k(u;u_1,...,u_{k-1})$. Then the family is defined by the quadratic
recursion relations:
\begin{eqnarray}
&\cQ_0(\bu)=1,\ \cQ_1(\bu) = u^{-1},& \nonumber \\ & \cQ_{k+1}(\bu)
= \frac{\displaystyle{\cQ_k^2(\bu)
-1}}{\displaystyle{u_k \cQ_{k-1}(\bu)}}, &\quad (k\geq
1). \label{quadratic}
\end{eqnarray}

Note that this system is the same as the $Q$-system if $u_i=1$ for
$i\geq 1$. Therefore, $\cQ_j(u,1,...,1) = Q_j$ where $u=t_1^{-1}$ to
agree with the initial conditions of \eqref{qsys}. Moreover, if we set
$u_1=\cdots = u_j=1$ leaving the other variables unevaluated, then the
solution of the system has $\cQ_k=Q_k$ for $k\leq j$ and $\cQ_{j+1} =
Q_{j+1}/u_j$.

The solutions to this system are known to be the Chebyshev polynomials
of the second kind.

\begin{lemma}\label{recursionlemma}
A family of solutions $\{\cQ_k\}$ satisfies \eqref{quadratic} if and
only if it satisfies the following recursion relations:
\begin{eqnarray}
& \cQ_0(\bu)=1,\ \cQ_1(\bu) = u^{-1},&
\nonumber \\ &\cQ_{k+1}(\bu) = \cQ_k(\bu'),&\ (k\geq 1),\label{recursion}
\end{eqnarray}
where 
$$
u' = \frac{1}{\cQ_2(\bu)}, \ u_1' = \cQ_1(\bu)u_2,\ u_j' = u_{j+1},\ (j>1),
$$
and $\cQ_2(\bu)$ is defined by the equation
\begin{equation}
\cQ_2(\bu) = \frac{u^{-2}-1}{u_1}.
\end{equation}
\end{lemma}

\begin{proof}
 Note that the quadratic relation in \eqref{quadratic} can be expressed as
\begin{equation}
\cQ_{k+1}(\bu) = \cQ_2(\cQ_k(\bu)^{-1},  \cQ_{k-1}(\bu)u_k).\label{Vtwo}
\end{equation}

Suppose the family of solutions satisfies the recursion
\eqref{recursion}. Then the quadratic relation \eqref{quadratic} holds
for $k=1$ by definition.
Suppose \eqref{quadratic} holds for $\cQ_m$ for all $m\leq k$, that is
$$
\cQ_m(\bu) = \cQ_2(\frac{1}{\cQ_{m-1}(\bu)},\cQ_{m-2}(\bu)u_{m-1}), \ m\leq k.
$$
Then
\begin{eqnarray*}
\cQ_{k+1}(\bu) &=& \cQ_k(\bu')\ \hbox{(by assumption)}\\
&=&
\cQ_2(\cQ_{k-1}(\bu')^{-1},\cQ_{k-2}(\bu')u_{k-1}') \ \hbox{(by induction hypothesis)}
\\
&=& \cQ_2(\cQ_k(\bu)^{-1},\cQ_{k-1}(\bu)u_k). 
\end{eqnarray*}
By induction, it follows that the family defined by \eqref{recursion}
satisfies the quadratic relation \eqref{quadratic} for all $k\geq 1$.

Conversely, suppose we have a family of functions which satisfies the
quadratic relation \eqref{quadratic}. Equation \eqref{recursion} holds
for $k=1$ by definition. 
Suppose \eqref{recursion} holds for all $m\leq k$. Then
\begin{eqnarray*}
\cQ_{k+1}(\bu) &=& \cQ_2(\cQ_k(\bu)^{-1},\cQ_{k-1}(\bu)u_k) \\
&=& \cQ_2(\cQ_{k-1}(\bu')^{-1},\cQ_{k-2}(\bu')u_{k-1}') \\
&=& \cQ_k(\bu').
\end{eqnarray*} 
The lemma follows by induction.
\end{proof}

\begin{remark}
The evaluation of the functions $\cQ_k$ when $u_j=1$ are the Chebyshev
polynomials of the second kind in the variable $t=u^{-1}$. These
polynomials, which are defined by the $sl_2$ fusion relation
\begin{eqnarray*}
& U_0(t) = 1, \ U_1(t) = t& \\
&U_1(t) U_k(t) = U_{k-1}(t) + U_{k+1}(t)& 
\end{eqnarray*}
are known to satisfy the quadratic relation 
$$
U_{k+1}(t) = \frac{U_k^2(t) - 1}{U_{k-1}(t)}.
$$ This is the $Q$-system for $\sl_2$, which is satisfied by the
characters of the irreducible representations of $sl_2$. Here,
$t=e^{\omega_1}+e^{-\omega_1}$.
\end{remark}

Lemma \ref{recursionlemma} can be recast in slightly more general terms.
Define the $j$th shift operation on the variables $\bu$ as:
\begin{equation}
u^{(j)} = \frac{1}{\cQ_{j+1}(\bu)}, \quad u_1^{(j)} = \cQ_j(\bu) u_{j+1},\quad
u_l^{(j)} = u_{l+j}\ (l>1)\label{uell}.
\end{equation}
The variable $\bu'$ is just $\bu^{(1)}$. Then
\begin{cor}
\begin{equation}\label{jtranslation}
\cQ_{k+j} (\bu) = \cQ_k(\bu^{(j)}).
\end{equation}
\end{cor}
\begin{proof}
By induction. The lemma holds for any $k$ when $j=1$ by Lemma
\ref{recursionlemma}. Suppose it is true for any $k$ and for all $l<j$. Then
\begin{eqnarray*}
\cQ_{k+j}(\bu) = \cQ_{k+j-1}(\bu') = \cQ_k((\bu')^{(j-1)})
\end{eqnarray*}
where the second equality is the induction hypothesis. We compute
\begin{eqnarray*}
(u')^{(j-1)} &=& \frac{1}{\cQ_j(\bu')} = \frac{1}{\cQ_{j+1}(\bu)} =
u^{(j)} \\
(u_1')^{(j-1)} &=& \cQ_{j-1}(\bu') u_{j+1}' = \cQ_j(\bu) u_{j+2} =
u_1^{(j)} \\
(u_m')^{(j-1)} &=& u_{m+1}^{(j-1)} = u_{m+j} = u_m^{(j)}.
\end{eqnarray*}
The statement follows.
\end{proof}

\subsection{Generating functions}
The general technique of the proof is to relax the restrictions on the
summations by defining an appropriate generating function. It is then
easy to prove properties of this generating function. In particular
this allows us to prove a constant term identity among generating
functions which is just the $M=N$ identity.

We define generating functions in the variables $\bu$, which are labeled
by the parameters $k\in \N$ and $\mathbf n\in \Z_+^k$ and $l\in \Z_+$:
\begin{equation}\label{genfun}
Z_{l;\mathbf n}^{(k)}(\mathbf u) := \sum_{m\in \Z_+^k} u^{q}
\prod_{i=1}^{k} {m_i + q_i \choose m_i} u_i^{q_i},
\end{equation}
where the integers $q_i$ depend on $m_i$ and $n_i$
$$q=l+\sum_{j=1}^k j(2 m_j-n_j), \quad q_i =q + p_i =l+ \sum_{j=1}^{k-i}
j(2m_{i+j} -n_{i+j}).$$ 
In particular, $q_k=l$.  Here, the binomial
coefficient is defined for all $p\in\Z$ by
$$
{m+p \choose m} = \frac{(p+m)(p+m-1) \cdots (p+1)}{m!}.
$$
Note that this coefficient is non-vanishing if
$m<-p$ in the case that $p$ is negative.  

This generating function is constructed so that the constant term in
$u$ corresponds to $q=0$. In this term, $q_i=p_i$, and the evaluation
at $u_i=1$ for all $i$ is just the HKOTY $N$-sum in equation
\eqref{mainidentity}.

On the other hand, by first taking terms with only non-negative powers
of each of variables $u_i$ and then taking the evaluation at $u_i=1$
for all $i$ and considering the constant term in $u$, we obtain the
$M$ side of equation \eqref{mainidentity}.

We will prove a stronger result, which is an identity of power series in $u$,
which implies the equality of these two types of constant terms. Then
in the limit $k\to\infty$, this will prove the combinatorial
Kirillov-Reshetikhin conjecture for $\g=\sl_2$.

The generating function $Z_{l;\mathbf n}^{(k)}(\mathbf u) $ is a
Laurent series in each $u_i$ and in $u$.  (The dependence on $u_k$ is
trivial, it is an overall factor $u_k^l$.)  Note that $q_j$ does not
depend on $m_i$ with $i\leq j$. It is therefore possible to find a
factorization formula for $Z_{l;\mathbf n}^{(k)}(\mathbf u) $ by
summing over $m_1$, then $m_2$ and so forth. In fact, such a
factorization can be described rather nicely. To prove it, we first
prove a simple Lemma.

\begin{lemma}
The function $ Z_{l;\mathbf n}^{(k)}(\mathbf u)$ satisfies the recursion 
relation
\begin{equation}
Z_{l;\bn}^{(k)}(\bu) = Z_{0;n_1}^{(1)}(\bu) Z_{l;n_2,...,n_k}^{(k-1)}(\bu')
\end{equation}
\end{lemma}
\begin{proof}
The function with $k=1$ can be computed explicitly (note that $q_1=l$
in this case)
\begin{eqnarray}
Z_{l,n_1}^{(1)}(u) &=& \sum_{m_1} u^{2m_1-n_1+l}u_1^{l}{m_1+l\choose
m_1}\nonumber\\ &=& \frac{u_1^l u^{-n_1+l}}{(1-u^2)^{l+1}}\nonumber \\
&=& \frac{\cQ_1(\bu)^{n_1+l+2}}{u_1 \cQ_2(\bu)^{l+1}}.\label{Zkone}
\end{eqnarray}

The summation over $m_1$ in \eqref{genfun} can be computed by using
the expansion (true for $|t|<1$ and for any $p$)
\begin{equation}\label{binomialsum}
\sum_{m\geq 0}{m+p\choose m} t^m = \frac{1}{(1-t)^{p+1}}.
\end{equation}
The result is
\begin{eqnarray*}
Z_{l,\bn}^{(k)}(\bu) &=& \sum_{m_2,...,m_k} u^{2q_1-q_2} u_1^{q_1}
\prod_{i=2}^{k} {m_i+q_i\choose m_i} u_i^{q_i}\times \sum_{m_1\geq 0}
u^{2m_1-n_1}{m_1+q_1\choose m_1}\\
&=& \frac{u^{-n_1}}{1-u^2} \sum_{m_2,...,m_k} \left(\frac{u^2
    u_1}{1-u^2}\right)^{q_1} u^{-q_2} \prod_{i=2}^{k}
{m_i+q_i\choose m_i} u_i^{q_i} \\
&=& Z_{0;n_1}^{(1)}(u) \sum_{m'_1,...,m_{k-1}'} {u'}^{q'}
\prod_{i=1}^{k-1} {m'_i+q'_i\choose m'_i} {u'_i}^{q'_i},
\end{eqnarray*}
where we have used the fact that 
$$2 q_1-q_2=2\sum_{j=2}^k j m_j = q - (2m_1-n_1).$$
Here, we have used the substitutions
$$
q'_i =l+ \sum_{j=1}^{k-i-1} j(2m_{j+i}'-n_{j+i}'), \quad n'_{j}=n_{j+1},\ m'_j = m_{j+1},
$$
and $q'=q_1$. The new variables $\bu'$ are
\begin{eqnarray*}
u' &=& \frac{u^2 u_1}{1-u^2} = \frac{1}{\cQ_2(\bu)} \\
u'_1 &=& \cQ_1(\bu) u_2 \\
u'_j &=& u_{j+1}, \ j>1.
\end{eqnarray*}
This proves the recursion relation for the generating function.
\end{proof}

This recursion allows us to prove the factorization formula for the
generating function $Z_{l;\bn}^{(k)}(\bu)$:
\begin{thm}\label{sltwoZfactorization}
The generating function $Z_{l;\bn}^{(k)}(\bu)$ has a factorization in
terms of the functions $\cQ_i(\bu)$
\begin{equation}\label{factorization}
Z_{l;\bn}^{(k)}(\bu) = \frac{\cQ_1(\bu)
  \cQ_k(\bu)^{l+1}}{\cQ_{k+1}(\bu)^{l+1}}\prod_{i=1}^k\frac{\cQ_i(\bu)^{n_i}
  }{ u_i} .
\end{equation}
\end{thm}
\begin{proof}
We prove the Theorem by induction. For $k=1$, the statement of the Theorem is
equivalent to equation \eqref{Zkone}. 
Suppose the Theorem is true for $k-1$. Then
\begin{eqnarray*}
Z_{l;\bn}^{(k)} = Z_{0,n_1}^{(1)}(\bu) Z_{l;n_2,...,n_k}^{(k-1)}(\bu')&=&
\frac{{\cQ_1}^{2+n_1}(\bu) }{\cQ_2(\bu) u_1} \times \frac{\cQ_1(\bu')
  \cQ_{k-1}(\bu')^{l+1}} {\cQ_{k}(\bu')^{l+1}} \prod_{i=1}^{k-1}\frac{
\cQ_i(\bu')^{n_{i+1}}}{ (u'_i)}\\
&=& \frac{{\cQ_1}^{2+n_1}(\bu) }{\cQ_2(\bu) u_1}\times \frac{\cQ_2(\bu)
  \cQ_k(\bu)^{l+1}}{\cQ_{k+1}(\bu)^{l+1}\cQ_1(\bu)} \prod_{i=2}^k \frac{
\cQ_i(\bu)^{n_i}}{u_i } \\
&=& \frac{\cQ_1(\bu) \cQ_k(\bu)^{l+1}}{\cQ_{k+1}(\bu)^{l+1}} 
\prod_{i=1}^k \frac{
\cQ_i(\bu)^{n_i}}{u_i},
\end{eqnarray*}
and the factorization is true for $k$. By induction, the factorization
holds for all $k\in\N$.
\end{proof}

\begin{cor}
Given any $1\leq p\leq k-1$, we have the factorization
\begin{equation}\label{pfactorization}
Z_{l;\bn}^{(k)}(\bu) = 
Z_{0;n_1,...,n_p}^{(p)}(\bu) Z_{l;n_{p+1},...,n_k}^{(k-p)}(\bu^{(p)}).
\end{equation}
\end{cor}
\begin{proof}
This follows from the factorization formula \eqref{factorization}
and the property \eqref{jtranslation}
\begin{eqnarray*}
\hspace{-.5in}&&
Z_{0;n_1,...,n_p}^{(p)}(\bu)Z_{l;n_{p+1},...,n_k}^{(k-p)}(\bu^{(p)})\\
&&\hspace{.5in}=\frac{\cQ_1(\bu) \cQ_p(\bu)}{\cQ_{p+1}(\bu)} \prod_{i=1}^p
\frac{\cQ_i(\bu)^{n_i}}{u_i} \times
\frac{\cQ_1(\bu^{(p)}) \cQ_{k-p}(\bu^{(p)})^{l+1}}{\cQ_{k-p+1}(\bu^{(p)})^{l+1}}
\prod_{i=1}^{k-p} \frac{ \cQ_i(\bu^{(p)})^{n_i}}{u_i^{(p)}} \\
&&\hspace{.5in}= \frac{\cQ_1(\bu) \cQ_p(\bu)}{\cQ_{p+1}(\bu)} \prod_{i=1}^p
\frac{\cQ_i(\bu)^{n_i}}{u_i}\times
\frac{\cQ_{p+1}(\bu)\cQ_{k}(\bu)^{l+1}}{\cQ_p(\bu)\cQ_{k+1}(\bu)^{l+1}}
\prod_{i=1}^{k-p}\frac{ \cQ_{i+p}(\bu)^{n_{i+p}}}{u_{i+p}} \\
&&\hspace{.5in}= \frac{\cQ_1(\bu) \cQ_{k}(\bu)^{l+1}}{\cQ_{k+1}(\bu)^{l+1}} \prod_{i=1}^k
\frac{\cQ_i(\bu)^{n_i}}{ u_i} = Z_{l;\bn}^{(k)}(\bu).
\end{eqnarray*}
\end{proof}

\subsection{Identity of power series}

Now consider the evaluation $\varphi_j$, which maps each variable in
the subset $\{u_1,...,u_{j-1}\}\subset \{u_1,...,u_k\},\ (j\leq k)$ to the
value $1$. That is,
$$\varphi_j (\bu) = (u;1,...,1,u_j,...,u_{k}).$$
At this specialization we have
\begin{eqnarray*}
\varphi_j ( \cQ_i(\bu)) &=& U_i(u^{-1}),\ i=1,...,j,\\
\varphi_j (\cQ_{j+1}(\bu)) &=& \frac{U_{j+1}(u^{-1})}{u_j}.
\end{eqnarray*}
Using the factorization \eqref{pfactorization} and the fact that
$$\varphi_j(Z_{0;n_1,...,n_j}^{(j)}(\bu)) =
\frac{U_1(u^{-1})U_j(u^{-1})}{U_{j+1}(u^{-1})}\
\prod_{i=1}^j U_i(u^{-1})^{n_i},
$$
we have
\begin{equation}\label{jfact}
\varphi_j(Z_{l;n_1,...,n_k}^{(k)}(\bu)) = 
\left(\frac{U_1(u^{-1})U_j(u^{-1})}{U_{j+1}(u^{-1})} 
\prod_{i=1}^j U_i(u^{-1})^{n_i}\right)\times \varphi_j(
Z_{l;n_{j+1},...,n_k}^{(k-j)} (\bu^{(j)})).
\end{equation}
Here, 
$$ \varphi_j(\bu^{(j)}) = (\frac{u_j}{U_{j+1}(u^{-1})}; U_j(u^{-1})
u_{j+1}, u_{j+2},...,u_k).
$$

Note that all the dependence on the parameters $u_j,...,u_{k}$ in
$\varphi_j(Z_{l;n_1,...,n_k}^{(k)}(\bu))$ 
is in the second factor of \eqref{jfact}. 

\begin{defn}
Let $f(u)$ be a Laurent series in $u$. Denote the power series part
(non-negative powers) of a Laurent series by
$$
\cP_u f(u).
$$
\end{defn}
\begin{defn}
Let $Z_{l;\bn}^{(k)}(\bu)$ be defined by equation \eqref{genfun}, and
let $1\leq j \leq k-1$. We denote by $Z_{l;\bn}^{(k)}(\bu)^{[j]}$ the
generating functions in \eqref{genfun} with the summation restricted
to values of $\bm$ such that $q_i\geq 0$ for all $i\geq j$. Note that
this is equivalent to taking only non-negative powers in the variables
$u_i$ with $i\geq j$.
\end{defn}

It is clear from the factorization formula that for any
$i\geq j$,
\begin{eqnarray}\label{evalfactorization}
\varphi_j ( Z_{l;n_1,...,n_k}^{(k)}(\bu)^{[i]})
=
\left[\frac{U_1(u^{-1})U_j(u^{-1})}{U_{j+1}(u^{-1})} 
\prod_{i=1}^j U_i(u^{-1})^{n_i}\right]
\ \varphi_j(
Z^{(k)}_{l;n_{j+1},...,n_k}(\bu^{(j)})^{[i]}).
\end{eqnarray}

\begin{lemma}\label{inductiveZ}
\begin{equation}\label{powerseries}
\cP_{u}
\varphi_j(Z_{l;n_1,...,n_k}^{(k)}(\bu)^{[j]}) =
\cP_{u}
\varphi_j(Z_{l;n_1,...,n_k}^{(k)}(\bu)^{[j+1]}).
\end{equation}
In other words, the power series in $u$ on the right hand side, where
the integers $q_{j+1},...,q_{k-1}$ are restricted to non-negative
values, has only non-negative powers of $u_{j}$.
\end{lemma}
\begin{proof}
First consider the factorization formula \eqref{jfact} with
$j=k-1$: 
\begin{eqnarray*}
\varphi_{k-1}(Z_{l;n_1,...,n_k}^{(k)}(\bu)) &=& 
\varphi_{k-1}(Z_{0;n_1,...,n_{k-1}}^{(k-1)}(\bu)) 
\varphi_{k-1}(Z_{l;n_k}^{(1)}(\bu^{(k-1)}))\\
&& \hskip-1in=
\frac{U_1(u^{-1}) U_{k-1}(u^{-1})}{U_k(u^{-1})}
\prod_{i=1}^{k-1}U_i(u^{-1})^{n_i} \
\sum_{m_k\geq 0} \left(\frac{u_{k-1}}{U_k(u^{-1})}\right)^{2m_k-n_k+l}
{m_k+l\choose m_k}u_k^l U_{k-1}(u^{-1})^l.
\end{eqnarray*}
We analyze the dependence on $u$ as follows. On the right hand side,
terms corresponding to strictly {\em negative} powers of $u_{k-1}$ are
proportional to a product of Chebyshev polynomials in $u^{-1}$, since
the factor $U_{k}$ in the denominator cancels. 

This means that the
coefficient of $u_{k-1}^{-n}\ (n>0)$ is a polynomial in $u^{-1}$.
Moreover, this polynomial has an overall factor of
$U_1(u^{-1})=u^{-1}$, and thus contains no constant term in $u^{-1}$.
Therefore, terms with strictly negative powers in $u_{k-1}$ appear
only in the coefficients of strictly negative powers of $u$ in the
Laurent series $\varphi_{k-1}(Z_{l;n_1,...,n_k}^{(k)}(\bu))$.

In terms of the notation above we have shown that
\begin{equation}\label{basestep}
\cP_u \varphi_{k-1} (Z_{l;n_1,...,n_k}^{(k)}(\bu))=
\cP_{u} \varphi_{k-1}(Z_{l;n_1,...,n_k}^{(k)}(\bu)^{[k-1]}).
\end{equation}

We proceed to prove the Lemma by induction. Consider the series
obtained from Equation \eqref{evalfactorization} by taking power series of both sides
\begin{eqnarray*}
&& 
\varphi_j(Z_{l;n_1,...,n_k}^{(k)}(\bu)^{[j+1]})= \\
&& \hskip.5in \frac{U_1(u^{-1})
  U_{j}(u^{-1})}{U_{j+1}(u^{-1})}\prod_{i=1}^j U_i(u^{-1})^{n_i} \times
\cP_{u_{j+1},...,u_{k-1}}\varphi_j(
Z_{l;n_{j+1},...,n_k}^{(k-j)}({\bu}^{(j)}))
\end{eqnarray*}
Explicitly,
\begin{eqnarray*}
& & \cP_{u_{j+1},...,u_{k-1}}\varphi_j(
Z_{l;n_{j+1},...,n_k}^{(k-j)}(\bu^{(j)}))\\
&&\hskip.5in=\sum_{\underset{q_{s}\geq 0 (s>j)}{m_{j+1},...,m_k\geq 0}}
\left(\frac{u_j}{U_{j+1}}\right)^{q_j} U_j^{q_{j+1}}\prod_{i=j+1}^{k}
  {q_i+m_i\choose m_i} u_i^{q_i}.
\end{eqnarray*}
Again, it is clear all terms with  strictly negative powers of $u_j$
($q_j<0$) are proportional to a product of Chebyshev polynomials in $u^{-1}$,
and this is where all the dependence on the variable $u$ resides. 
Because of the overall factor $U_1(u^{-1})$, such terms contribute
only to the coefficients of $u^{-n}$ with $n>0$ in the generating function. 

The Lemma follows by induction.
\end{proof}

\begin{thm}\label{sltwoHKOTY}
\begin{equation}\label{HKOTY}
\cP_u Z_{l;\bn}^{(k)}(u,1,...,1)^{[1]}  = 
\cP_u Z_{l;\bn}^{(k)} (u,1,...,1).
\end{equation}
\end{thm}
\begin{proof}
Lemma \ref{inductiveZ} guarantees that 
\begin{equation}\label{indstep}
 \cP_u Z_{l;\bn}^{(k)}(u,1,...,1))^{[j]}=\cP_u
Z_{l;\bn}^{(k)}(u,1,...,1)^{[j+1]},
\end{equation} 
by evaluating both sides of \eqref{powerseries} at the point
$u_{j}=...=u_{k}=1$. 

We proceed by induction, with the initial step
coming from equation \eqref{basestep}:
$$
\cP_u Z_{l;\bn}^{(k)}(u,1,...,1)) = \cP_u Z_{l;\bn}^{(k)}(u,1,...,1))^{[k-1]}.
$$
The induction step is \eqref{indstep}.
The Theorem follows.
\end{proof}

The relation to the restriction in the summation of the HKOTY
conjecture is as follows. The constant term in the generating function
$Z_{l,\bn}^{(k)}$ corresponds to all terms with $q=0$. This constant
term appears as the first term of the power series identity in $u$
proven above.

Keeping in mind that $q_i = p_i + q = p_i$ in the constant term,
the constant term of the right hand side of
\eqref{HKOTY} is the unrestricted ($N$)-side of the HKOTY conjecture,
and the constant term on the left hand side of \eqref{HKOTY} is the
restricted ($M$)-side of the HKOTY conjecture. Thus, we have proven the
conjecture for the case $\g=\sl_2$.

\subsection{Identity of multiplicities}

The number $N_{l;\bn}$ is the multiplicity of the irreducible
$\sl_2$-module with highest weight $l\omega_1$ in the tensor product
$\otimes_i V(i \omega_1)^{n_i}$. This number is equal to the
multiplicty of the trivial representation in $V(l \omega_1) \otimes
\left( \otimes_i V(i \omega_i)^{n_i}\right)$. 
As a non-trivial check, we show that the above generating function
indeed gives this property.

The factorized form of $Z_{l;\bn}^{(k)}(\bu)$ allows for writing a
simple residue integral for the multiplicity $N_{l;\bn}$ of the
representation $V_l$ in the tensor product $\otimes_{i=1}^p
V_i^{\otimes \, n_i}$.
\begin{lemma}\label{multiplicitysltwo}
The multiplicities $N_{l;\bn}$ are equal to the residue integral around $u=0$:
\begin{equation}
N_{l;\bn}=\oint \frac{du}{2i\pi u} \prod_{i=1}^p U_i(u^{-1})^{n_i} \, U_1(u^{-1}) z(u)^{l+1}\label{multitwo}
\end{equation}
where $z(u)=uC(u^2)$, $C(x)=(1-\sqrt{1-4x})/(2x)=\sum_{n\geq 0}c_n x^n$ being the generating series of the 
Catalan numbers $c_n=(2n)!/(n!(n+1)!)$.
\end{lemma}
\begin{proof}
The integer $N_{l;\bn}$ is the constant term of
$Z_{l;\bn}^{(k)}(u,1,1,...,1)$ in the limit as $k\to \infty$, while
only finitely many $n_i$ are nonzero. Assume $\bn$ is such that
$n_j=0$ for all $j>p$ for some $p\in \N$ (we pick $p$ so that $l\leq p$).
  Then
\begin{eqnarray*}
N_{l;\bn}&=&\lim_{k\to \infty} \oint \frac{du}{2i\pi u} Z_{l;\bn}^{(k)}(u,1,1,...,1)\\
&=&\oint {du \over 2i\pi u} \prod_{i=1}^p U_i(u^{-1})^{n_i} \lim_{k\to \infty}
U_1(u^{-1})\frac{U_k(u^{-1})^{l+1}}{U_{k+1}(u^{-1})^{l+1}}
\end{eqnarray*}
where the contour of integration encircles $0$.
With the parametrization $u^{-1}=z+z^{-1}$,
where $z=z(u):=(2u)^{-1}(1-\sqrt{1-4u^2})$, the Chebyshev polynomials read 
$U_k(u^{-1})=(z^{k+1}-z^{-k-1})/(z-z^{-1})$, and we have the following large $k$ 
asymptotics for $|z|<1$: $U_k(u^{-1})\sim -z^{-k}/(z^2-1)$, so that 
$\lim_{k\to \infty} U_k(u^{-1})/U_{k+1}(u^{-1}) =z(u)$. The lemma follows, as $z(u)/u=C(u^2)$.
\end{proof}

As a non-trivial check of the formula \eqref{multitwo}, we show that
$N_{l;\bn}=N_{0;\bn+\epsilon_l}$, where the vector
${\bf \epsilon}_l$ has components $\delta_{i,l}$ for $i\geq 1$.

Using the expression
\eqref{multitwo} for $N_{0;\bn+{\bf \epsilon}_l}$ and the fact that
$U_l (u^{-1})=z^l +z^{-1} U_{l-1}(u^{-1})$, we may rewrite
\begin{eqnarray*}
N_{0;\bn+{\bf \epsilon}_l}&=&\oint \frac{du}{2i\pi u} \prod_{i=1}^p U_i(u^{-1})^{n_i} \, U_1(u^{-1})U_l(u^{-1}) z(u) \\
&=& \oint \frac{du}{2i\pi u} \prod_{i=1}^p U_i(u^{-1})^{n_i} \, U_1(u^{-1})\Big( z(u)^{l+1}+ U_{l-1}(u^{-1})\Big)\\
&=& N_{l;\bn}
\end{eqnarray*}
as the second term $U_{l-1}(u^{-1})$ yields a product of polynomials
of $u^{-1}$ with $U_1(u^{-1})$ as an overall factor, hence has no
constant term in $u$, while the first term exactly matches
eq. \eqref{multitwo} for $N_{l;\bn}$.

\section{The HKOTY identity for $\g$ a simply-laced Lie algebra}

The treatment illustrated for $\g=\sl_2$ in the previous section is
essentially unchanged for other simple Lie algebras, although with
more complicated notation. It is simplest to describe the non
simply-laced algebras separately. Thus, we give the proof of
\eqref{mainidentity} in this section for the case of simply-laced Lie
algebras. The arguments of the previous section generalize in a
straightforward manner to these algebras.

Let $\g$ be a simple, simply-laced Lie algebra with Cartan matrix $C$,
and rank $r$. 

\subsection{The $\cQ$-system}
Let $\bu$ denote the set of formal variables $\bu:=\{u_\al;
u_{\al,i}|\ 1\leq \al\leq r, i\in \N\}$. 

The family of functions
$\{\cQ_{\alpha,k}(\bu): \ 1\leq \alpha\leq r, k\in \Z_+\}$ is defined
recursively as follows:
\begin{eqnarray}\label{gquadratic}
&& \cQ_{\al,0}(\bu)=1,\quad \cQ_{\al,1}(\bu)=(u_{\al})^{-1},\nonumber \\ 
&&\cQ_{\al,k+1}(\bu) =
  \frac{{(\cQ_{\al,k}(\bu))^2 - \prod_{\beta\neq\alpha} 
    (\cQ_{\beta,k}(\bu))^{-C_{\beta,\alpha}}} }{\displaystyle{u_{\al,k} \cQ_{\al,k-1}(\bu)}}.
\end{eqnarray}

\begin{remark}\label{reduction}
Note that if $u_{\al,j}=1$ for all $\al$ and all $j$, then the
$\cQ$-system is equivalent to the $Q$-system for simply-laced $\g$,
with given initial conditions for $Q_{\al,1}$. The functions of
$u_{\al}^{-1}$ thus defined are ``generalized Chebyshev polynomials''
in the variables $u_{\al}^{-1}.$
\end{remark}

In particular,
\begin{equation}\label{vtwo}
\cQ_{\al,2}(\bu) =\frac{(1-\prod_{\beta}
(u_{\beta})^{-C_{\beta,\al}})}{u_{\al}^2 u_{\al,1}}. 
\end{equation}
Therefore, the quadratic relation \eqref{gquadratic} can be expressed as
\begin{equation}
\cQ_{\al,j+1}(\bu) = \cQ_{\al,2}(\bu^{(j-1)}),
\end{equation}
where
\begin{equation}\label{utranslate}
u_{\al}^{(j)} = \frac{1}{\cQ_{\al,j+1}(\bu)},\ 
u_{\al,1}^{(j)} = \cQ_{\al,j}(\bu) u_{\al,j+1},\ 
u_{\al,l}^{(j)} = u_{\al,l+j},\ (l>1).
\end{equation}

\begin{lemma}
A family of functions $\{\cQ_{\al,k}: \ 1\leq \alpha\leq r, k\in
\Z_+\}$ satisfies \eqref{gquadratic} if and only if it satisfies the
following recursion relation:
\begin{eqnarray}\label{grecursion}
\cQ_{\al,0}=1,\ \cQ_{\al,1}=u_{\al}^{-1},\nonumber\\
\cQ_{\al,k+1}(\bu) =
\cQ_{\al,k}(\bu'),
\end{eqnarray}
where $\bu' = \bu^{(1)}$ is defined by equation \eqref{utranslate} and $\cQ_2(\bu)$ is defined by equation \eqref{vtwo}.
\end{lemma}
\begin{proof}
Suppose the family of functions $\{\cQ_{\al,k}\}$ satisfies
\eqref{gquadratic}. Then $\cQ_{\al,2}(\bu)$ satisfies \eqref{grecursion}
because
$$
\cQ_{\al,2}(\bu) = \frac{1}{u_{\al}'}= \cQ_{\al,1}(\bu') 
$$
from the definition \eqref{utranslate} for $j=1$. 
Suppose \eqref{grecursion} holds for all $\cQ_{\al,m}(\bu)$ with
$m\leq k$. Then
\begin{eqnarray*}
\cQ_{\al,k+1} &=& \cQ_{\al,2}(\bu^{(k-1)})\qquad \hbox{(by
  \eqref{gquadratic})} 
\\ 
&=& \cQ_{\al,2}((\bu')^{(k-2)})\qquad \hbox{(by induction hypothesis)} \\
&=& \cQ_{\al,k}(\bu') \qquad \hbox{(by \eqref{gquadratic})}.
\end{eqnarray*}
Here, we used the fact that $(\bu')^{(k-2)} = \bu^{(k-1)}$, which
follows from the induction hypothesis, because
$$
(u_{\al}')^{(k-2)} = \frac{1}{\cQ_{\al,k-1}(\bu')} = \frac{1}{\cQ_{\al,k}(\bu)}
$$
and so forth.
By induction, the recursion \eqref{grecursion} holds for all $k$.

Conversely, suppose the family $\{\cQ_{\al,k}(\bu)\}$ satisfies
\eqref{grecursion}. Then 
again the relation \eqref{grecursion} holds for $\cQ_{\al,2}(\bu)$ by
definition. Suppose it holds for all $m\leq k$. Then
\begin{eqnarray*}
\cQ_{\al,k+1}(\bu) &=& \cQ_{\al,k}(\bu') \qquad \hbox{(by \eqref{grecursion})}\\
&=& \cQ_{\al,2} (\bu'^{(k-2)})\qquad \hbox{(by induction hypothesis)} \\
&=& \cQ_{\al,2}(\bu^{(k-1)}),
\end{eqnarray*}
so that \eqref{gquadratic} holds for $\cQ_{\al,k+1}(\bu)$. The lemma
follows by induction.
\end{proof}

\begin{cor}
\begin{equation}\label{genrecursion}
\cQ_{\al,k+j}(\bu) = \cQ_{\al,k}(\bu^{(j)}).
\end{equation}
\end{cor}
\begin{proof}
The case $j=1$ is the statement of the last Lemma. Suppose
\eqref{genrecursion} holds for a fixed $k$ for all $l<
j$. Then 
$$
\cQ_{\al,k+j}(\bu) = \cQ_{\al,k+j-1}(\bu') =
\cQ_{\al,k}((\bu')^{(j-1)}) = \cQ_{\al,k}(\bu^{(j)}). 
$$
and \eqref{genrecursion} holds for all $j$ by induction.
\end{proof}

In the specialization to the case where $u_{\al,i}=1$ for $i>0$, the
$\cQ$-system degenerates to the $Q$-system for $\g$
simply-laced. Therefore, the specialization of $\cQ_{\al,j}(\bu)$ to
this point gives $Q_{\al,j}$.
A theorem of Nakajima \cite{Nakajima} (see
also \cite{Hernandez}) shows that the solution to the $Q$-system is
the set of characters of the KR-modules of $\g$. That is, $Q_{\al,j}$
is the character of $\KR_{\al,j}$.


\subsection{Generating functions}
Fix $k\in\N$. 
Given a set of non-negative integers $\{m_{\al,i},n_{\al,i}:\ 1\leq
i\leq k, 1\leq \alpha\leq r\}$ and a $\g$-dominant integral weight $\lambda =
\sum_\al l_\al \omega_\al$, , define for each $\al\in I_r$ and $i\in \N$
\begin{eqnarray}\label{gspin}
q_{\al} &=& l_\al + \sum_{j=1}^k \sum_{\beta\in I_r} j
(C_{\al,\beta}m_{\beta,j} - \delta_{\al,\beta} n_{\beta,j})\label{totalspin}\\
q_{\al,i}&=& q_\al + p_{\al,i} = l_\al +\sum_{j=1}^{k-i} \sum_{\beta\in I_r} j(
C_{\alpha,\beta}m_{\beta,i+j}-\delta_{\alpha,\beta} n_{\beta,i+j}).\label{vacancy}
\end{eqnarray}
Below, we will use the notation $\bm =
(m_{1,1},m_{1,2}...,m_{1,k},m_{2,1},...)$, $\bm_1=\{m_{1,1},...,m_{r,1}\}$ and so forth. Note that if
we extend equation \eqref{vacancy} to $i=0$ then $q_\al=q_{\al,0}$. In
the case that $i=k$ there is no dependence on the parameters $\bn,\bm$
and $q_{\al,k}=l_\al$.

Define the generating function in the variables $\{u_\al;u_{\al,i}|\
1\leq i\leq k, 1\leq \alpha\leq r\}$:
\begin{equation}\label{gZ}
Z_{\lambda;\bn}^{(k)}(\bu) = 
\sum_{\bm \in Z_+^{r\times k}} \prod_{\al=1}^r u_\al^{q_\al}
\prod_{i=1}^k {m_{\al,i}+q_{\al,i}\choose m_{\al,i}} u_{\al,i}^{q_{\al,i}}.
\end{equation}
This generating function is constructed so that the constant term (with
$q_\al=0$ for all $\al$) is the $N$ side \eqref{N} of the identity \eqref{mainidentity} when
all $u_{\al,i}=1$, because in the constant term, $q_{\al,i}=p_{\al,i}$.

Note that, for any fixed set of values of $q_\al, (\al\in I_r)$, the
coefficient of $\prod_\al u_\al^{q_\al}$ is a Laurent polynomial in the
variables $\{u_{\al,i}\}$. The generating series is thus 
a Laurent series in the variables $\{u_{\al}\}$, with coefficients which are Laurent
polynomials in the other the variables.
In fact, because of the form \eqref{gspin} of $q_\al$, the dependence on
the variables $u_\al$ is such that, up to an overall factor depending on
$\bn$ and $\lambda$, $Z_{\lambda,\bn}^{(k)}(\bu)$ is a power series
in the variables 
$$
y_\al = \prod_\beta u_{\al}^{C_{\al,\beta}},
$$
since
$$\prod_\al \prod_\beta u_\al^{C_{\al,\beta}\sum_j jm_{\beta,j}} =
\prod_\beta \left[\prod_\al u_\al^{C_{\al,\beta}}\right]^{\sum_j
 jm_{\beta,j}} = \prod_\beta y_\beta^{\sum_j jm_{\beta,j}}.$$

\begin{lemma}
\begin{equation}\label{Zgone}
Z_{\lambda;\bn_1}^{(1)}(\bu) = \prod_{\al=1}^r
\frac{u_{\al,1}^{l_\al}u_{\alpha}^{-n_{\alpha,1}+l_\al}}
{(1-\prod_{\beta=1}^r u_{\beta}^{C_{\alpha,\beta}})^{l_\al+1}}. 
\end{equation}
where we interpret the denominator of the term corresponding to $\al$ on the right hand side
as a power series in the variable
$$
y_\al=\prod_\beta u_{\beta}^{C_{\al,\beta}}, \quad (\al\in I_r).
$$
\end{lemma}
Note that expanding as a power series in $y_\al$ is equivalent to expanding in $u_\al$ for each $\al$ since $y_\al$ has non-negative powers only in $u_\al$.
\begin{proof}
We use the definition \eqref{gZ}, noting that $q_{\al,1}=l_\al$ in
this case:
$$ Z_{\lambda;\bn_1}^{(1)} = \sum_{m_{1,1},...,m_{r,1}\geq 0}
\prod_{\al} u_\al^{q_\al} { m_{\al,1}+l_\al\choose
m_{\al,1}} u_{\al,1}^{l_\al}.
$$
Here, $q_{\al,1} = l_\al + \sum_\beta m_{\beta,1} - n_{\al,1}$. Thus,
\begin{eqnarray*}
Z_{\lambda;\bn_1}^{(1)} = \prod_\al
u_\al^{-n_{\al,1}+l_\al}u_{\al,1}^{l_\al} \sum_{m_{1,1},...,m_{r,1}}
\prod_\al \left( \prod_\beta
u_\beta^{C_{\al,\beta}}\right)^{m_{\al,1}}{m_{\al,1}+l_\al\choose
m_{\al,1}}.
\end{eqnarray*}
We can perform the summation over each $m_{\al,1}$ using equation
\eqref{binomialsum}. The Lemma follows.
\end{proof}
Note that 
$Z_{\lambda;\bn_1}^{(1)}(\bu)$ 
can be written in terms of the functions $\cQ(\bu)$:
\begin{equation}\label{initialZ}
Z_{\lambda;\bn_1}^{(1)}(\bu) = \prod_{\al=1}^r 
\frac{\cQ_{\al,1}(\bu)^{n_{\al,1}+l_\al+2}}
{u_{\al,1} \cQ_{\al,2}(\bu)^{l_\al+1}}.
\end{equation}
Again, in each factor,  the denominator corresponding to $\al$ is to be regarded as a power series
expansion in the variable $u_\al$.

We now proceed as in the $\sl_2$ case.
\begin{lemma}\label{zgrecursion}
The function $Z_\bn^{(k)}(\bu)$ satisfies the following recursion relation:
\begin{equation}
Z_{\lambda;\bn}^{(k)}(\bu) = 
Z_{0;\bn_1}^{(1)}(\bu) Z_{\lambda;\bn_2,...,\bn_k}^{(k-1)}(\bu'),
\end{equation}
where $\bu'$ is defined by equation \eqref{utranslate} with $j=1$.
\end{lemma}
\begin{proof}
Note the identity
\begin{equation}\label{qgreaterthanone}
q_\al^{(1)} := q_{\al}\Big|_{\bm_1=\bn_1=0} = 2
  q_{\al,1}-q_{\al,2},
\end{equation}
Therefore,
\begin{eqnarray*}
Z_{\lambda;\bn}^{(k)}(\bu) &=& \sum_{\bm_2,...,\bm_k}
\prod_{\al}u_{\al,1}^{q_{\al,1}}
u_{\al}^{q_{\al}^{(1)}}\prod_{i=2}^k u_{\al,i}^{q_{\al,i}}{ q_{\al,i}+m_{\al,i}\choose
  m_{\al,i}} \sum_{\bm_1} \prod_\al
u_{\al}^{\sum_\beta C_{\al,\beta}m_{\beta,1} -
n_{\alpha,1}}{q_{\al,1}+m_{\al,1}\choose m_{\al,1}},
\end{eqnarray*}
Since the functions $\{q_{\al,1}\}_{\al\in I_r}$ do not depend on $\bm_1$, the sum over
$\bm_1$
can be performed explicitly:
\begin{eqnarray*}
Z_{\lambda;\bn}^{(k)}(\bu) &=&\prod_\al\frac{
u_{\al}^{-n_{\al,1}}}{1-\prod_\beta u_{\beta}^{C_{\al,\beta}}}
\sum_{\bm_2,...,\bm_k}
\prod_{\al}\left[\frac{u_{\al}^2u_{\al,1}}{1-\prod_\beta
u_{\beta}^{C_{\al,\beta}}}\right]^{q_{\al,1}}
u_{\al}^{-q_{\al,2}}\prod_{i=2}^k u_{\al,i}^{q_{\al,i}}{
q_{\al,i}+m_{\al,i}\choose m_{\al,i}}\\
&=& Z_{0;\bn_1}^{(1)}(\bu) \sum_{\bm_2,...,\bm_k}
\prod_{\al}\left[\frac{1}{\cQ_{\al,2}}\right]^{q_{\al,1}}\cQ_{\al,1}^{q_{\al,2}} \prod_{i>1}
u_{\al,i}^{q_{\al,i}}{
q_{\al,i}+m_{\al,i}\choose m_{\al,i}}.
\end{eqnarray*}
We have used equation
\eqref{Zgone} with $\lambda=0$ to identify the first factor.  
The second factor is the
generating function $Z_{\lambda,\bn'}^{(k-1)}(\bu')$, where
$\bn'_i=\bn_{i+1}$ with $i=1,...,k-1$ and $\bu'$ are the variables
defined via the substitution \eqref{utranslate} with $j=1$.
\end{proof}

The main factorization theorem is:
\begin{thm}\label{gfactorization}
There is a factorization
\begin{equation}
Z_{\lambda;\bn}^{(k)}(\bu) = 
\prod_{\al=1}^r \cQ_{\al,1}(\bu) 
\left(\frac{\cQ_{\al,k}(\bu)}{\cQ_{\al,k+1}(\bu)}\right)^{l_\al+1}
\prod_{i=1}^k \frac{\cQ_{\al,i}(\bu)^{n_{\al,i}}}{u_{\al,i}}.
\end{equation}
Here, each $\cQ_{\al,i}$ in the denominator is understood as a Laurent series expansion in $u_\al$ for each $\al$.
\end{thm}
\begin{proof}
The Lemma is true for $k=1$ by equation \eqref{initialZ}. Assume it is
true for $k-1$. Use the recursion of Lemma \ref{zgrecursion}:
\begin{eqnarray*}
Z_{\lambda;\bn}^{(k)}(\bu) &=& Z_{0;\bn_1}^{(1)}(\bu)
Z_{\lambda;\bn'}^{(k-1)} (\bu') \\
&=& \prod_{\al=1}^r \frac{\cQ_{\al,1}(\bu)^{n_{\al,1}+2}}
{u_{\al,1} \cQ_{\al,2}(\bu)} \times
\cQ_{\al,1}(\bu')
\left[\frac{\cQ_{\al,k-1}(\bu')}{\cQ_{\al,k}(\bu')}\right]^{l_\al+1}
\prod_{i=2}^k \frac{\cQ_{\al,i-1}(\bu')^{n_{\al,i}}}{u'_{\al,i-1}} \\
&=&\prod_\al \frac{\cQ_{\al,1}(\bu)^{n_\al + 2}\cQ_{\al,2}(\bu)}
{u_{\al,1} \cQ_{\al,2}(\bu)} \left[
    \frac{\cQ_{\al,k}(\bu)}{\cQ_{\al,k+1}(\bu)}\right]^{l_\al+1}
  \frac{1}{\cQ_{\al,1}(\bu)} \prod_{i=2}^k
  \frac{\cQ_{\al,i}(\bu)^{n_\al,i}}{u_{\al,i}} \\
&=& \prod_\al 
  \cQ_{\al,1}\left[\frac{\cQ_{\al,k}(\bu)}{\cQ_{\al,k+1}(\bu)}
\right]^{l_\al+1} \prod_{i=1}^k
  \frac{\cQ_{\al,i}(\bu)^{n_{\al,i}}}{u_{\al,i}}, 
\end{eqnarray*}
where we have used the fact that
$\cQ_{\al,i}(\bu')=\cQ_{\al,i+1}(\bu)$.
\end{proof}

\begin{cor}
There is a factorization formula
\begin{equation}\label{gfact}
Z_{\lambda;\bn}^{(k)}(\bu) = 
Z_{0;\bn_1,...,\bn_j}^{(j)}(\bu) Z_{\bn^{(j)}}^{(k-j)}
(\bu^{(j)}),
\end{equation}
where $\bu^{(j)}$ is defined by equation \eqref{utranslate} and $\bn^{(j)}=(\bn_{j+1},...,\bn_k)$.
\end{cor}
\begin{proof}
The proof uses the factorization of Theorem \ref{gfactorization}.
\begin{eqnarray*}
Z_{\lambda;\bn}^{(k)}(\bu) &=& \prod_{\alpha} \frac{\cQ_{\al,1}(\bu)
  \cQ_{\al,k}(\bu)^{l_\al+1}}{\cQ_{\al,k+1}(\bu)^{l_\al+1}} \prod_{i=1}^j
\frac{\cQ_{\al,i}(\bu)^{n_{\al,i}}}{u_{\al,i}} \times \prod_{i=j+1}^k 
\frac{\cQ_{\al,i}(\bu)^{n_{\al,i}}}{u_{\al,i}}\\
&&\hskip-.5in= \prod_{\alpha} \left[
  \left(\frac{\cQ_{\al,1}(\bu)\cQ_{\al,j}(\bu)}{\cQ_{\al,j+1}(\bu)}
\prod_{i=1}^j\frac{\cQ_{\al,i}(\bu)^{n_{\al,i}}}{u_{\al,i}}\right)\left(
\frac{\cQ_{\al,k}(\bu)^{l_\al+1}\cQ_{\al,j}(\bu)}{\cQ_{\al,k+1}(\bu)^{l_\al+1}}\prod_{i=1}^{k-j}
\frac{\cQ_{\al,i}(\bu^{(j)})^{n_{\al,i+j}}}{u_{\al,i}^{(j)}}\right)\right] \\
&&\hskip-.5in= \prod_{\alpha} \left[
  \left(\frac{\cQ_{\al,1}(\bu)\cQ_{\al,j}(\bu)}{\cQ_{\al,j+1}(\bu)}
\prod_{i=1}^j\frac{\cQ_{\al,i}(\bu)^{n_{\al,i}}}{u_{\al,i}}\right)\left(
\frac{\cQ_{\al,k-j}(\bu^{(j)})^{l_\al+1}\cQ_{\al,1}(\bu^{(j)})}
{\cQ_{\al,k-j+1}(\bu^{(j)})^{l_\al+1}}\prod_{i=1}^{k-j}
\frac{\cQ_{\al,i}(\bu^{(j)})^{n_{\al,i+j}}}{u_{\al,i}^{(j)}}\right)\right]\\
&&= Z_{0;\bn_1,...,\bn_{j}}(\bu) Z_{\lambda;\bn_{j+1},...,\bn_{k}}(\bu^{(j)}).
\end{eqnarray*}
\end{proof}

\subsection{Identity of power series}

We need a Lemma about series expansions of $(Q_{\al,j})^{-1}$
in the variables $y_\al=u_\al^2 \prod_{\beta\sim \al} u_\beta^{-1}$.
\begin{lemma}\label{seriesform}
Let $A=\Z[u_1^{-1},...,u_{r}^{-1}]$. If we
interpret $1/Q_{\al,2}$ as an element of the ring $A[[u_\al]]$,
then for all $j\geq 2$, 
$$
\frac{1}{Q_{\al,j}}\in A[[u_\al]].
$$
\end{lemma}
\begin{proof}
This statement follows from the particular form of the $Q$-system and from the fact that $Q_{\al,j}\in A $ by Theorem \ref{polynomiality}.
Note that for $j=2$, we interpret
$$
\frac{1}{Q_{\al,2}} = u_\al^2 (1-y_\al)^{-1}= u_\al^2\sum_{m\geq 0}y_\al^m  \in u_\al^2 A[[y_\al]] = A[[u_\al]].
$$ Then for any $j> 2$, negative
powers of $Q_{\al,j}$ are in the same ring.
For suppose  $Q_{\al,j}^{-m}\in A[[u_\al]]$. Then using the $Q$-system,
$$
\frac{1}{Q_{\al,j+1}} =
\frac{Q_{\al,j-1}}{Q_{\al,j}^2(1-\prod_\beta Q_{\beta,j}^{-C_{\al,\beta}})}.
$$
By the induction hypothesis and Theorem \ref{polynomiality}, 
$$
\frac{Q_{\al,j-1}}{Q_{\al,j}^2}\in A[[u_\al]].
$$
Moreover, 
$$ \prod_\beta Q_{\beta,j}^{-C_{\al,\beta}} =
Q_{\al,j}^{-2}\prod_{\beta\sim\al} Q_{\beta,j}\in A[[u_\al]]
$$
for the same reason. Therefore,
the denominator is expanded as a function in $A[[u_\al]]$, that is,
$$
\frac{1}{1-\prod_\al Q_{\al,j}^{-C_{\al,\beta}}} = \sum_{m\geq 0}
\left(\prod_\al Q_{\al,j}^{-C_{\al,\beta}}\right)^{m}.
$$

The Lemma follows by induction.
\end{proof}

\begin{defn}
The evaluation map $\varphi_j$ acting on functions of $\bu$ is defined by
$$ \varphi_j(u_{\al,l})=\left\{ \begin{array}{ll} 1,& l<j \\ u_{\al,l}
& \hbox{otherwise},\end{array}\right.
$$ extended by linearity.
\end{defn}

We use the fact, which follows from the definition
\eqref{gquadratic} and remark \ref{reduction}, that
\begin{equation}
\varphi_j(\cQ_{\alpha,l}(\bu)) = Q_{\al,l}, (l\leq j),\qquad
\varphi_j(\cQ_{\al,j+1}(\bu))=\frac{Q_{\al,j+1}}{u_{\al,j}}.
\end{equation}
This follows from the fact that
$\varphi_j(\cQ_{\al,2}(\bu)),...,\varphi_j(\cQ_{\al,j}(\bu))$ satisfy the usual
$Q$-system in this case, and the dependence of $\varphi_j(\cQ_{\al,j+1}(\bu))$ on
$u_{\al,j}$ is explicitly just an overall factor, and it otherwise
also satisfies the usual $Q$-system.

\begin{defn}
Let $Z_{\lambda;\bn}^{(k)}(\bu)$ be defined as the summation over
$\bm$ as in equation \eqref{gZ}. Let $K=(k_1,...,k_r)\in \N^r$. Define
$Z_{\lambda;\bn}^{(k)}(\bu)^{[K]}$ to be the restricted sum generating function,
defined as in \eqref{gZ} but  with the summation over $\bm$ restricted to
values of $q_{\al,j}\geq 0$ for all $j\geq k_\al$ for each $\al$. 
\end{defn}
Recall that $\cP_{u_\al}$ denotes the power series part of the Laurent expansion in the variable $u_\al$ of a function of $u_\al$.

\begin{lemma}\label{powerseriesg} Fix a root $\al$.
Let $K=(k_1,...,k_r)$  with $k_\beta\geq j$ for $\beta\neq \al$  and $k_\al=j+1$. Define $\epsilon_\al$ to be an $r$-vector with $1$ in the $\al$th entry. Then
\begin{equation}\label{powerserieseq} 
\cP_{u_\al} 
\varphi_j(Z_{\lambda;\bn}^{(k)}(\bu)^{[K]}) = 
\cP_{u_\al} 
\varphi_j(Z_{\lambda;\bn}^{(k)}(\bu)^{[K-\epsilon_\al]}).
\end{equation}
\end{lemma}
\begin{proof}
Consider the factorization formula \eqref{gfact}.  The restriction
$[K]$ does not affect the first factor, so the factorization formula
still holds for the restricted summation.  We prove the Lemma by
induction. The base step is to take $j=k-1$. Let
$K=(k_1,...,k_\al,...,k_r)$ with $k_\beta\geq k-1$ for each $\beta$,
and $k_\al=k$
\begin{eqnarray}\label{base}
\varphi_{k-1} (Z_{\lambda;\bn}^{(k)}(\bu)^{[K]}) &=& \varphi_{k-1} (Z_{0;\bn_1,...,\bn_{k-1}}^{(k-1)}(\bu)) \varphi_{k-1}
(Z_{\lambda;\bn_k}^{(1)}(\bu^{(k-1)}) ^{[K]})\nonumber\\ 
&=&
\left[\prod_\al\frac{Q_{\al,1} Q_{\al,k-1}}{Q_{\al,k}} \prod_{i=1}^{k-1}
Q_{\al,i}^{n_{\al,i}}\right] \varphi_{k-1}(
Z_{\lambda;\bn_k}^{(1)}(\bu^{(k-1)})^{[K]}) \\
&=& \prod_\al\frac{Q_{\al,1} Q_{\al,k-1}}{Q_{\al,k}} \prod_{i=1}^{k-1}
Q_{\al,i}^{n_{\al,i}} \nonumber\\
&&\times\sum_{\bm_k\atop q_{\beta,i}\geq 0 \ (i\geq k_\beta)}
\prod_\al \left[\frac{u_{\al,k-1}}{Q_{\al,k}}\right]^{q_{\al,k-1}}
{m_{\al,k}+l_{\al}\choose m_{\al,k}} (Q_{\al,k-1}u_{\al,k})^{l_\al}.\nonumber
\end{eqnarray}
Note that we have assumed that $l_\al\geq 0$.

Consider a term in the summation over $\bm_k$, with fixed
$q_{\al,k-1}$ such that $ q_{\al,k-1}<0$.  The dependence on the
variables $u_\al$ in equation \eqref{base} is only via the functions
$\{Q_{\beta,i}\}_{\beta\in I_r, i\in \N}$.  That is, in the factor
$$ \prod_{\beta\in I_r}
Q_{\beta,k-1}^{1+l_\beta}Q_{\beta,1}Q_{\beta,k}^{-q_{\beta,k-1}-1}\prod_{i=1}^{k-1}
Q_{\beta,i}^{n_{\beta,i}}.$$ Moreover, terms with non-negative powers
in $u_\al$ can only come from the expansion in the factor
corresponding to $\beta=\alpha$. If $q_{\al,k-1}<0$ this is a
polynomial in the $Q_{\al,i}$ for various $i$, which is in
$A$. Moreover, it comes with a positive overall power of
$Q_{\al,1}=u_{\al}^{-1}$ and hence has no constant term in $u_{\al}$
at all. Therefore the power series expansion in $u_\al$ has no
contribution from this term.

We have shown that
$$
\cP_{u_\al} \varphi_{k-1}(Z_{\lambda;\bn}^{(k)}(\bu))=
\cP_{u_\al}\varphi_{k-1}(Z_{\lambda;\bn}^{(k)}(\bu)^{[K]})=
\cP_{u_\al} \varphi_{k-1}(Z_{\lambda;\bn}^{(k)}(\bu)^{[K-\epsilon_\al]}).
$$

Next, suppose the Lemma is true for $j+1$. Note that a restriction on
the integers $q_{\beta,i}$ with $i\geq j$ does not involve the
summation over the integers $m_{\al,i}$ with $i\leq j$. Therefore
there is still a partial factorization of the generating function with
restrictions into a product of two factors. That is, if
$K=(k_1,...,k_r)$ with $k_\beta \geq j$ and $k_\al = j+1$, then
$$ Z_{\lambda;\bn}^{(k)}(\bu)^{[K]} = Z_{0;\bn_1,...,\bn_j}^{(j)}
(\bu)Z_{\lambda;\bn^{(j)}}^{(k-j)}(\bu^{(j)})^{[K]}.
$$ One should read the restriction in the second factor as a
restriction to non-negative powers of $u_{\beta,i}$ with $i\geq
k_\beta$.

Therefore, we can write
\begin{eqnarray*}
\varphi_j (
Z_{\lambda,\bn}^{(k)}(\bu)^{[K]}) &=& 
\varphi_j(
Z_{0;\bn_1,...,\bn_j}^{(j)}(\bu)) \varphi_j(
Z_{\lambda;\bn^{(j)}}^{(k-j)}(\bu^{(j)}) ^{[K]})\\ 
&=&
\prod_\al\frac{Q_{\al,1} Q_{\al,j}}{Q_{\al,j+1}} \prod_{i=1}^j
Q_{\al,i}^{n_{\al,i}} \times \varphi_j(Z_{\lambda;\bn^{(j)}}^{(k-j)}(\bu^{(j)})^{[K]}) .
\end{eqnarray*}
We use the definition for the second factor:
$$
\varphi_j(
Z_{\lambda;\bn^{(j)}}^{(k-j)}(\bu^{(j)})^{[K]})
=
\sum_{\underset{q_{\beta,s}\geq 0 \ (s\geq k_\beta)}{\bm^{(j)}}}
\prod_\al
\left[\frac{u_{\al, j}}{Q_{\al,j+1}}\right]^{q_{\al,j}} 
Q_{\al,j}^{q_{\al,j+1}}
\prod_{i=j+1}^k {m_{\al,i}+q_{\al,i}\choose m_{\al,i}} u_{\al,i}^{q_{\al,i}}.
$$

Consider a term in the summation with $q_{\al,j}<0$ and
$q_{\al,j+1}\geq 0$ for some $\al$. The dependence on $u_\al$ is
contained in the functions $Q_{\beta,i}$. Moreover, non-negative
powers of $u_\al$ can only come from the expansion of $Q_{\al,j+1}$
appearing in the denominator. If $q_{\al,j}<0$ there are no such
terms, and we are left with a polynomial in the $Q_{\al,i}$s with an
overall factor $Q_{\al,1}$, so there is no constant term in $u_\al$,
and we have a polynomial in $u_{\al}^{-1}$.

We have
$$
\cP_{u_\al} \varphi_j (Z_{\lambda;\bn^{(j)}}^{(k-j)}(\bu^{(j)})^{[K]})
=\cP_{u_\al} \varphi_j (Z_{\lambda;\bn^{(j)}}^{(k-j)}(\bu^{(j)})^{[K-\epsilon_\al]})
$$
where $K$ has $k_\beta\geq j$ and $k_\al=j+1$.
The Lemma follows by induction on $j$.

\end{proof}
We have the obvious corollary:
\begin{cor}\label{alphaseries}
Let $J=\{\al_1,...,\al_t\}\subset I_r$. Let $K$ be a set with 
$k_{\al}=j+1$ for $\al\in J$, and $k_\al=j$ for $\al\in I_r\setminus J$.
Let $K'$ be the set with $k_\al=j$ for all $\al\in I_r$. Then we have
\begin{equation}\label{seriesalpha}
\cP_{u_{\al_1},...,u_{\al_t}} \varphi_{j}
(Z_{\lambda;\bn}^{(k)}(\bu)^{[K]}) = 
\cP_{u_{\al_1},...,u_{\al_t}} \varphi_{j} (Z_{\lambda;\bn}^{(k)}(\bu)^{[K']}).
\end{equation}
\end{cor}
\begin{proof}
This follows by repeating the argument in the Lemma above for several
values of $\al$.
\end{proof}

This implies an identity of power series:
\begin{thm}\label{gHKOTY}
\begin{equation}
\cP_{u_1,...,u_r}
\varphi_{k+1}(Z_{\lambda;\bn}^{(k)}(\bu)) = 
 \cP_{u_1,...,u_r}\varphi_{k+1}(
Z_{\lambda;\bn}^{(k)}(\bu)^{[1,...,1]}). 
\end{equation}
\end{thm}
\begin{proof}
This follows by induction on $j$ from the previous corollary, setting $t=r$.
The evaluation map $\varphi_{k+1}$ can then be applied to both sides of the
resulting power series identity.
\end{proof}
\begin{cor} The conjecture \ref{HKOTYconj} is true for $\g$ simply laced.
\end{cor}
\begin{proof}
Theorem \ref{gHKOTY} is an identity of power series in
$\{u_\al\}_{\al\in I_r}$. The constant term of this identity (the
restriction to $q_\al=0$ for all $\al$) is the $M=N$ conjecture of
\cite{HKOTY} for the case of $\g$ simply-laced.
\end{proof}
\section{The identity for non-simply laced Lie algebras}

\subsection{The functions $\cQ_{\alpha,k}$}
Let $\g$ be one of the non simply-laced algebras.  It is convenient to
define the integer $t_\al$, which takes the value 1 if $\al$ is a long
root, $t_\al=2$ if $\al$ is a short root of $B_r, C_r$ or $F_4$ and
$t_\al=3$ for the short root ($\al=2$) of $G_2$.

We define a family of functions $\{\cQ_{\alpha,i}(\bu)|\ 
(\al,i)\in I_r\times \Z_+\}$, depending on the formal variables
$\bu=\{u_\al, \ u_{\al,j}, a_j\ | \ \al\in I_r, \ j\in \N \}$,
via the recursion
relations:
\begin{eqnarray}\label{quadraticQ}
&&\cQ_{\al,0}(\bu)=1,\quad \cQ_{\al,1}(\bu)=u_\al^{-1}, \nonumber \\
&&\cQ_{\al,i+1}(\bu) = \frac{\cQ_{\al,i}(\bu)^2 -\prod_{\beta\sim\al} \mathcal
T_i^{(\al,\beta)}(\bu)} {u_{\al,i} \cQ_{\al,i-1}(\bu)},\ i\geq1,
\end{eqnarray}
with
\begin{equation*}
 \mathcal T_i^{(\al,\beta)}(\bu) =
a_{i}^{-i({\rm mod}t_\al)}
\prod_{k=0}^{|C_{\al,\beta}|-1}
\cQ_{\beta,\lfloor{\frac{t_\beta i+k}{t_\al}\rfloor}}(\bu).
\end{equation*}

Explicitly, we have
$$
\mathcal T_i^{(\al,\beta)}(\bu) = \cQ_{\beta,i}(\bu)\quad
\hbox{if $t_\al=t_\beta$}
$$ 
and
\begin{eqnarray}
& B_r: & \cT_i^{(r-1,r)} = \cQ_{r,2i},\label{exceptions} \\
& & \cT_{2i-1}^{(r,r-1)} = a_{2i-1} \cQ_{r-1,i-1}\cQ_{r-1,i},\nonumber\\
& & \cT_{2i}^{(r,r-1)} = \cQ_{r-1,i}^2.\nonumber\\
& C_r: & \cT_{2i-1}^{(r-1,r)} = a_{2i-1}\cQ_{r,i-1}\cQ_{r,i},\nonumber\\
&& \cT_{2i}^{(r-1,r)} = \cQ_{r,i}^2,\nonumber\\
&& \cT_{i}^{(r,r-1)} = \cQ_{r-1,2i}.\nonumber\\
& F_4: & \cT_{i}^{(2,3)} = \cQ_{3,2i},\nonumber\\
&& \cT_{2i-1}^{(3,2)} = a_{2i-1} \cQ_{2,i-1}\cQ_{2,i},\nonumber\\
&& \cT_{2i}^{(3,2)} =  \cQ_{2,i}^2.\nonumber\\
& G_2: & \cT_i^{(1,2)} = \cQ_{2,3i},\nonumber\\
&& \cT_{3i-2}^{(2,1)} = a_{3i-2}^2 \cQ_{1,i} \cQ_{1,i-1}^2,\nonumber\\
&& \cT_{3i-1}^{(2,1)} = a_{3i-1} \cQ_{1,i}^2 \cQ_{1,i-1}.\nonumber\\
&& \cT_{3i}^{(2,1)} =  \cQ_{1,i}^3,\nonumber
\end{eqnarray}

Note that when $a_i=u_{\al,i}=1$ for all $i$ and $\al$, this is just
the $Q$-system of Section \ref{q-system}, with the identification of
the initial conditions. Thus, we have a deformed $Q$-system, which we
refer to as the $\cQ$-system from now on.

\begin{lemma}\label{usefulindeed}
A family of functions $\{\cQ_{\al,i} | \ \al\in I_r, i\in \Z_+ \}$
satisfies the $\cQ$-system \eqref{quadraticQ} if and only if it
satisfies the recursion relation
\begin{equation}\label{recQ}
\cQ_{\al,i+t_\al}(\bu) = \cQ_{\al,i}(\bu'), \ i\geq 1,
\end{equation}
subject to the initial conditions that $\cQ_{\al,j} (j\leq t_\al)$ are
defined by \eqref{quadraticQ}, and such that
\begin{equation}\label{shiftedua}
u_{\al}' = \frac{1}{\cQ_{\al,t_\al+1}},\quad u_{\al,j}' =
\cQ_{\al,t_\al}^{\delta_{j,1}}u_{\al,t_\al+j},\quad
a_{i}' = \cQ_{\gamma,1}^{\delta_{i<t_{\gamma'}}}a_{i+t_{\gamma'}}
\end{equation}
Here, $\gamma=r-1$ for $B_r$, $\gamma=r$ for $C_r$, 
$\gamma = 2$ for $F_4$ and $\gamma = 1$ for $G_2$, and $\gamma'$ is
the short root connected to $\gamma$, i.e. $\gamma'=r,r-1,3,2$ for
$\g=B_r,C_r,F_4,G_2$, respectively. The function $\delta_{i<j}$ is 1
if $i<j$ and 0 otherwise.
\end{lemma}

\begin{proof}

Depending on the value of $t_\al$, the proof has up to three steps,
but proceeds in a similar way as for the simply-laced case.

Suppose that a family of functions $\{\cQ_{\al,j}(\bu)\}$ satisfies the
$\cQ$-system \eqref{quadraticQ} for all $(\al,j)$. 
We will show that it also satisfies the recursion relation \eqref{recQ}.

The initial conditions of both systems are the same by definition.
Therefore, suppose that the functions $\cQ_{\al,j}(\bu)$ also satisfy the
recursion relation \eqref{recQ} for any $j\leq t_\al k$ for some
$k$. We will show 
that it also satisfies \eqref{recQ} for $j\leq t_\al(k+1)$.

For all $\g$, the function $\cT_{t_\al k}^{(\al,\beta)}(\bu')$ is some
power of the function
$\cQ_{\beta,t_\beta k}(\bu')$. The induction
hypothesis is that $\cQ_{\beta,t_\beta k}(\bu')
= \cQ_{\beta,t_\beta(k+1)}(\bu)$ Therefore,
\begin{equation}\label{recursionT}
\cT_{t_\al k}^{(\al,\beta)}(\bu') =
\cT_{t_\al(k+1)}^{(\al,\beta)}(\bu).
\end{equation}

Now, consider
\begin{eqnarray}
\cQ_{\al, t_\al k+1}(\bu') &=& \frac{\cQ_{\al,t_\al k }(\bu') -
\prod_{\beta\sim \al}
\cT^{(\al,\beta)}_{t_\al k}(\bu')}{u_{\al,t_\al k}' \cQ_{\al,t_\al k-1}
(\bu')} \quad \hbox{(by assumption, true for all $k$)}\nonumber \\ 
&=&
\frac{\cQ_{\al,t_\al(k+1)}(\bu) - \prod_{\beta\sim \al}
\cT^{(\al,\beta)}_{t_\al k}(\bu')} {u_{\al,t_\al (k+1)} 
\cQ_{\al,t_\al (k+1)-1}(\bu)}\quad \hbox{(by induction
  hypothesis.)}\nonumber \\
&=&\cQ_{\al,t_\al(k+1)+1}(\bu)\qquad  \hbox{(by equation \eqref{recursionT}.)}
\label{plusone}
\end{eqnarray}

In the cases where $t_\al>1$, we must now consider the functions
$\cT_{t_\al k + 1}^{(\al,\beta)}(\bu')$.  If $\al=\gamma'$, this function is
proportional to some power of $a_{t_\al k+1}'=a_{t_\al (k+1)+1}$ 
(for $k\geq 1$). It
is also proportional to (a product of some powers of)
$\cQ_{\beta, j}(\bu')$ with $j \in \{t_\beta k,t_\beta k +1\}$.  By the
induction hypothesis and by \eqref{plusone}, we have that
$\cQ_{\beta,j}(\bu') = \cQ_{\beta,j+t_\beta}(\bu)$ for these values of $j$.
Therefore, we have that 
\begin{equation}\label{recursionTtwo}
\cT_{k t_\al+1}^{(\al,\beta)}(\bu') = \cT_{(k+1)
t_\al+1}^{(\al,\beta)}(\bu).
\end{equation}

We can now use \eqref{recursionTtwo} and \eqref{plusone}, together
with the induction hypothesis, to obtain, for $t_\al>1$,
\begin{eqnarray*}
\cQ_{\al,t_\al k + 2}(\bu') &=& \frac{\cQ_{\al,t_\al k + 1}(\bu') -
\prod_{\beta\sim \al} \cT_{t_\al k + 1}^{(\al,\beta)}(\bu')} {
  u_{\al,t_\al k + 1}' \cQ_{\al,t_\al k}(\bu')} \\
&=& \frac{\cQ_{\al,t_\al (k+1) + 1}(\bu) - \prod_{\beta\sim \al}
    \cT_{t_\al (k+1)+1}(\bu)}{u_{\al,t_\al(k+1) + 1}
    \cQ_{\al,t_\al(k+1)}(\bu)} \\
&=& \cQ_{\al,t_\al(k+1)+2}(\bu).
\end{eqnarray*}

Finally, consider the case of $t_\al=3$ for $\g=G_2$. In this case, 
$$
\cT_{3k+2}^{(2,1)} (\bu') = a_{3k+2}' \cQ_{1,k+1}(\bu')^2 \cQ_{1,k}(\bu')
=
a_{3k+5} \cQ_{1,k+2}(\bu)^2 \cQ_{1,k+1}(\bu) =
\cT_{3(k+1)+2}^{(2,1)}(\bu).
$$
It follows that for $G_2$, $\cQ_{2,3k+2}(\bu') =
\cQ_{2,3(k+1)+2}(\bu)$. We have thus proven that the family of
functions which satisfies the $\cQ$-system \eqref{quadraticQ} also
satisfies the recursion \eqref{recQ}.

Conversely, suppose we have a family of functions defined by the
recursion relations \eqref{recQ}, with the functions
$\cQ_{\al,j}(\bu)$ being identical to those satisfying the
$\cQ$-system for all $j\leq t_\al$. Then we show that the family
$\cQ_{\al,j}$ satisfies the $\cQ$-system for all $j$. 

Above, we showed
that if the functions $\cQ_{\al,j}$ satisfy the recursion relation
\eqref{recQ}, then the functions $\cT_j^{(\al,\beta)}$ satisfy a
similar recursion relation:
$$
\cT_j^{(\al,\beta)}(\bu') = \cT_{j+t_\al}^{(\al,\beta)}(\bu).
$$ Therefore, if $\cQ_{\al,j}$ satisfies the $\cQ$-system for all
$j\leq t_\al k$, then
\begin{eqnarray*}
\cQ_{\al,j+t_\al}(\bu) &=& \cQ_{\al,j}(\bu') \\ &=& \frac{
\cQ_{\al,j-1}(\bu')-\prod_{\beta\sim \al}
\cT_{j-1}^{(\al,\beta)}(\bu')} {u_{\al,j-1}' \cQ_{\al,j-2}(\bu')} 
\\ &=& \frac{\cQ_{\al,j+t_\al-1}(\bu) - \prod_{\beta\sim \al}
\cT_{j+t_\al-1}^{(\al,\beta)}(\bu)}{ u_{\al,j-1+t_\al}
\cQ_{\al,j+t_\al-2} (\bu)}.
\end{eqnarray*}
But this is just the $\cQ$-system for $\cQ_{j+t_\al}$.
Therefore, by induction, the $\cQ$-system is satisfied all $j, \al$.
\end{proof}

\begin{cor}\label{usefulcor}
The family of functions $\cQ_{\al,j}$ defined by either equation
\eqref{quadraticQ} or \eqref{recQ} satisfies the recursion
\begin{equation}
\cQ_{\al,k+t_\al j} (\bu)= \cQ_{\al,k}(\bu^{(j)}),
\end{equation}
where 
\begin{eqnarray}\label{uj}
u_{\al}^{(j)} &=& \frac{1}{\cQ_{\al,jt_\al + 1}},\quad
u_{\al,1}^{(j)}= \cQ_{\al, jt_\al} u_{\al,jt_\al+1},\nonumber \\
u_{\al,i}^{(j)} &=& u_{\al, i+ jt_\al}\ (i>1),\quad
a_{i}^{(j)}= \cQ_{\gamma,j}^{\delta_{i<t_{\gamma'}}} a_{i+j t_{\gamma'}}.
\end{eqnarray}
\end{cor}
\begin{proof}
By induction. The case $j=1$ is the case of Lemma
\ref{usefulindeed}. Suppose that the corollary holds for all $k$ and all
$j<l$. Then by Lemma \ref{usefulindeed}, we have
We note that $(\bu')^{(j)} = \bu^{(j+1)}$ from the definition, by
using Lemma \ref{usefulindeed}. Thus,
$$
\cQ_{\al,k+t_\al l}(\bu) = \cQ_{\al,k+t_\al(l-1)}(\bu')=
\cQ_{\al,k}((\bu')^{(l-1)}) = \cQ_{\al,k}(\bu^{(l)}).
$$
By induction on $l$, the Corollary follows.
\end{proof}

\subsection{Generating functions}
For a given natural number $k$, define the set of pairs $J_\g^{(k)} =
\{(\al,j): \al\in I_r, 1\leq j\leq t_\al k \}$.  Given a set of
non-negative integers $\{ m_{\al,i},n_{\al,i} \ | \ (\al,i)\in J_\g^{(k)}
\}$ and a dominant integral weight $\lambda=\sum_{\al=1}^r l_\al
\omega_\al$, $l_\al \in 
\Z_+$, we define the total spin as before:
$$
q_{\al} = l_\al+\sum_{(\beta, j)\in J_\g^{(k)} } j C_{\al,\beta} m_{\beta,j} 
-\sum_{j=1}^{t_\al k}j n_{\al,j}.
$$

In this case, the vacancy numbers have the form, for $(\al,i)\in
J_\g^{(k)}$:
$$ p_{\al,i} = \sum_{j=1}^{t_\al k} \min(i,j) n_{\al,j} - \sum_{(\beta,
  j)\in J_\g^{(k)}} {\rm 
sgn}(C_{\al,\beta}) \min(|C_{\al,\beta}|j,
|C_{\beta,\al}|i)m_{\beta,j}.
$$
For $(\al,i)\in J_\g^{(k)}$, define the modified vacancy numbers as for the
simply-laced case:
\begin{eqnarray}
 q_{\al,i} &=& p_{\al,i} + q_{\al} \nonumber \\ &=&
l_\al+\sum_{\underset{|C_{\al,\beta}|j>|C_{\beta,\al}|i}{(\beta, j) \in
J_\g^{(k)}:} } {\rm sgn}(C_{\alpha,\beta})
(|C_{\al,\beta}|j-|C_{\beta,\al}|i) m_{\beta,j} - \sum_{j=i+1}^{t_\al
k} (j-i) n_{\al,j}.\label{nonsimpleq}
\end{eqnarray}
Note that for any $\al$, $q_{\al,t_\al k}=l_\al.$ We list the explicit
forms of the modified vacancy numbers for each of the algebras in the
Appendix.

We define generating functions $Z_{\lambda;\bn}^{(k)}(\bu)$ parametrized
by sets of non-negative integers $\mathbf n := \{n_{\al,j}\ | \ (\al,j)\in
J_\g^{(k)}\}$ and a dominant integral weight $\lambda$ as follows:
\begin{equation}\label{zgen}
Z_{\lambda;\bn}^{(k)}(\bu) = \sum_{\bm}
\prod_{\al=1}^r u_\al^{q_\al} \prod_{i=1}^{t_\al k}
u_{\al,i}^{q_{\al,i}} a_i^{\Delta_{\al,i}}
{m_{\al,i}+q_{\al,i}\choose m_{\al,i}} ,
\end{equation}
where $\bm = \{m_{\al,i}\geq 0\ | (\al,i)\in J_\g^{(k)}\}$ and $\Delta_{\al,i}=(-i\ {\rm mod}\
t_{\gamma'})\delta_{\al,\gamma'}m_{\al,i}$. Explicitly, the
non-vanishing values of $\Delta_{\al,i}$ are
\begin{eqnarray*}
&B_r: & \Delta_{r,2i-1}=m_{r,2i-1} \\
&C_r: & \Delta_{r-1,2i-1}=m_{r-1,2i-1} \\
&F_4: & \Delta_{3,2i-1}=m_{3,2i-1} \\
&G_2: & \Delta_{2,3i-2}=2m_{2,3i-2},\quad \Delta_{2,3i-1}=m_{2,3i-1}.
\end{eqnarray*}

By a slight abuse of notation, throughout this section, we shall
always denote by $Z_{\lambda;\bn}^{(j)}(\bu)$ the sum \eqref{zgen}
with an appropriately truncated sequence $\bn$, namely, with
$n_{\al,i}=0$ for all $i>t_{\al} j$.

\subsubsection{Factorization properties of the generating function}

\begin{lemma}\label{lemmageninit}
For $k=1$, the generating function has the factorized form in terms of
the $\cQ$-functions:
\begin{equation}\label{initgen}
Z_{\lambda;\bn}^{(1)}(\bu) =\prod_{\al=1}^r
\frac{\cQ_{\al,1}\cQ_{\al,t_\al}^{l_\al+1}}
{\cQ_{\al,t_{\al}+1}^{l_\al+1}} \prod_{i=1}^{t_\al}
\frac{\cQ_{\al,i}^{n_{\al,i}}}{u_{\al,i}}
\end{equation}
where 
$\cQ_{\al,i}$ are defined by the $\cQ$-system \eqref{quadraticQ}. Here, the equality is understood as an equality of Laurent series in $u_\al$ (for each $\al\in I_r$), where each factor in the product over $\al$ should be expanded in the corresponding $u_\al$.
\end{lemma}
The proof is by direct calculation, which is done in section \ref{appthree} of the Appendix.

\begin{lemma}\label{recuZ}
The generating function $Z_{\lambda;\bn}^{(k)}$ satisfies the
recursion relation, for $k\geq 2$:
\begin{equation}\label{Zrecursion}
Z_{\lambda;\bn}^{(k)}(\bu) = Z_{0;\bn}^{(1)}(\bu)
Z_{\lambda;\bn'}^{(k-1)}(\bu'),
\end{equation}
where $\bu'$ is defined by equation \eqref{shiftedua}, and $\bn'$ is
obtained from $\bn$ omitting the integers $n_{\al,j}$ with $j\leq
t_\al$, and $u^{(j)}$ is defined by equation \eqref{uj}.
\end{lemma}

\begin{proof} 
By direct calculation. See Section \ref{lemmarecuZ} of the Appendix.
\end{proof}

\begin{thm}\label{Zfactorized}
The generating function $Z_{\lambda;\bn}^{(k)}(\bu)$ factorizes as
follows:
\begin{equation}\label{factogen}
Z_{\lambda;\bn}^{(k)} = \prod_{\al=1}^r
\frac{\cQ_{\al,1}\cQ_{\al,t_\al k}^{l_\al+1}}
{\cQ_{\al,t_{\al}k+1}^{l_\al+1}} \prod_{i=1}^{t_\al k}
\frac{\cQ_{\al,i}^{n_{\al,i}}}{u_{\al,i}},
\end{equation}
Where again, the equality is understood as an equality of Laurent series in $u_\al$ (for each $\al\in I_r$), where each factor in the product over $\al$ should be expanded in the corresponding $u_\al$.
\end{thm}

\begin{proof}
We proceed by induction. For $k=1$, the formula holds by Lemma
\ref{lemmageninit}.  Assuming it holds for $k-1$, we apply the
recursion hypothesis to eq.\eqref{Zrecursion} and then use Lemma
\ref{usefulindeed} to rewrite it as a function of $\bu$.
\end{proof}

\begin{cor} The generating function has a factorization as follows:
\begin{equation}\label{factojgen}
Z_{\lambda;\bn}^{(k)}(\bu) =Z_{0;\bn}^{(j)}(\bu)\, 
Z_{\lambda;\bn^{(j)}}^{(k-j)}(\bu^{(j)})
\end{equation}
where $\bn^{(j)}$ is obtained from $\bn$ by omitting the components 
$n_{\beta,i} \ (i\leq t_\beta j)$.
\end{cor}
\begin{proof} 
Using \eqref{factogen}, evaluate the ratio
\begin{equation} 
\frac{Z_{\lambda;\bn}^{(k)}(\bu)}{Z_{0;\bn}^{(j)}(\bu)}
=\prod_{\al=1}^r \frac{\cQ_{\al,t_\al j+1}\cQ_{\al,t_\al k}^{l_\al+1}}
{\cQ_{\al,t_{\al}k+1}^{l_\al+1}} 
\frac{\cQ_{\al,t_\al j+1}^{n_{\al,t_\al j+1}}}{\cQ_{t_\al j}u_{\al,t_\al j+1}}
\prod_{i=t_\al j+2}^{t_\al k} \frac{\cQ_{\al,i}^{n_{\al,i}}}{u_{\al,i}},
\end{equation}
then apply Corollary \ref{usefulcor} to express the r.h.s. as a function of $\bu^{(j)}$, easily identified with $Z_{\lambda;\bn^{(j)}}^{(k-j)}(\bu^{(j)})$ via \eqref{factogen}.
\end{proof}

We also need to consider the following generating
functions arising from ``partial summations'' 
$Z^{(k,p)}_{\lambda;\bn}$, for $0\leq p < \max\{t_\al\}$.
Let us denote by $\Pi^<$ be set of the
short simple roots of $\g$. Let $\Pi^>=\Pi\setminus \Pi^<$. 

\begin{defn}
For each $\g$, let 
$$J_\g^{(k,p)} = 
\{(\al,i)\ | \ 1\leq i\leq k\ (\al\in \Pi^>),\
p<i\leq t_\al k\ (\al\in \Pi^<)\}.$$
Define
\begin{eqnarray}
Z_{\lambda;\bn}^{(k,p)}(\bu) &=& \sum_{m_{\al,i}\geq 0\ (\al,i)\in
J_\g^{(k,p)}} \left[\prod_{\al\in \Pi^>} u_\al^{\tilde q_\al} \prod_{i=1}^k
u_{\al,i}^{q_{\al,i}}{q_{\al,i}+m_{\al,i}\choose
m_{\al,i}}\right]\nonumber\\ & & \hskip.5in\times
\left[ \prod_{\beta\in \Pi^<}
u_{\beta,p}^{q_{\beta,p}}\prod_{i=p+1}^{t_\beta k} 
a_{i}^{\Delta_{\beta,i}}
u_{\beta,i}^{q_{\beta,i}} {q_{\beta,i}+m_{\beta,i}\choose
  m_{\beta,i}}\right]
.
\end{eqnarray}
Here, $\tilde q_{\al} = q_{\al}|_{m_{\beta,i}=0\ (\beta,i)\notin
  J_{\g}^{(k,p)}}$.
\end{defn}
Note that $Z_{\lambda,\bn}^{(k,0)} = Z_{\lambda,\bn}^{(k)}$. Also note
that $Z_{\lambda;\bn}^{(k,p)}(\bu)$ does not depend on the entries
$n_{\al,i}$ with $i\leq p$ and $\al$ is a short root.

In terms of this partially summed generating function, we have a factorization
\begin{lemma}\label{partialfactorization}
\begin{equation}
Z^{(k)}_{\lambda;\bn}(\bu)=\prod_{\beta\in \Pi^<}
{\cQ_{\beta,1}\cQ_{\beta,p}\over \cQ_{\beta,p+1}}
\prod_{i=1}^p {\cQ_{\beta,i}^{n_{\beta,i}}\over u_{\beta,i}} 
Z^{(k,p)}_{\lambda;\bn}(\bu^{(0,p)})
\end{equation}
where
\begin{equation}\label{subszerop}
u^{(0,p)}_{\beta,p}={1\over \cQ_{\beta,p+1}},\
u^{(0,p)}_{\beta,p+1}=\cQ_{\beta,p}u_{\beta,p+1} \quad \forall\beta\in \Pi^<.
\end{equation}
The other components of $\bu$ are unchanged under the substitution.
\end{lemma}
\begin{proof} By direct calculation. See Section \ref{apptwo} in the Appendix.
\end{proof}

It is helpful to introduce, for fixed $\g$, $0\leq p<
\max\{t_\al\}$ and $j$, the notation 
\begin{equation}\label{deftau}
\tau_\al = \left\{\begin{array}{ll} t_\al j + p,& \al\in \Pi^<\\
j & \al \in \Pi^>.\end{array}\right.
\end{equation}

\begin{cor}\label{lastZfactor}
\begin{eqnarray}
Z_{\lambda;\bn}^{(k)}(\bu) &=& \left[\prod_{\al \in
I_r}{\cQ_{\al,1}\cQ_{\al,\tau_\al}\over \cQ_{\al,\tau_\al+1}}
\prod_{i=1}^{\tau_\al} {\cQ_{\al,i}^{n_{\al,i}}\over u_{\al,i}}\right]
Z_{\lambda;\bn^{(j,p)}}^{(k-j,p)}(\bu^{(j,p)})
\end{eqnarray}
where $\bu^{(j,p)} = (\bu^{(j)})^{(0,p)}$. That is,
\begin{equation}
u^{(j,p)}_{\beta, p}={1\over \cQ_{\beta,t_\beta j+p+1}},\ 
u^{(j,p)}_{\beta, p+1}=\cQ_{\beta,t_\beta
  j+p}u_{\beta,t_\beta j+p+1} \quad \forall \beta\in \Pi^< \  
\end{equation}
and the other components of $\bu^{(j,p)}$ are equal to those of
$\bu^{(j)}$. The set $\bn^{(j,p)}$ is the set $\bn$ without the
entries $n_{\al,i}$ with $i\leq \tau_\al$.
\end{cor}
\begin{proof}
From  the factorization \eqref{factojgen} and Lemma
\ref{partialfactorization}, we have 
\begin{eqnarray*}
Z_{\lambda;\bn}^{(k)} &=& Z_{0;\bn}^{(j)}(\bu)
Z_{\lambda;\bn^{(j)}}^{(k-j)}(\bu^{(j)})\\
&=&Z_{0;\bn}^{(j)}(\bu)\, \prod_{\beta\in \Pi^<}
{\cQ_{\beta,1}(\bu^{(j)})
\cQ_{\beta,p}(\bu^{(j)})\over \cQ_{\beta,p+1}(\bu^{(j)})}
\prod_{i=1}^p {\cQ_{\beta,i}^{n_{\beta,t_\beta j+i}}(\bu^{(j)})
\over u_{\beta,i}^{(j)}} \, 
Z_{\lambda;\bn^{(j,p)}}^{(k-j,p)}((\bu^{(j)})^{(0,p)})
\\ &=&
Z_{0;\bn}^{(j)}(\bu)\, 
\prod_{\beta\in \Pi^<}
{\cQ_{\beta,t_\beta j+1}\cQ_{\beta,t_\beta j+p}\over \cQ_{\beta,t_\beta j}
\cQ_{\beta,t_\beta j+p+1}}
\prod_{i=1}^p {\cQ_{\beta,t_\beta j+i}^{n_{\beta,t_\beta j+i}}\over 
u_{\beta,t_\beta j+i}} \, 
Z_{\lambda;\bn^{(j,p)}}^{(k-j,p)}(\bu^{(j,p)})
\end{eqnarray*}
where we have used Corollary \ref{usefulcor}. The corollary follows
from the factorization of Theorem \ref{Zfactorized} applied to the
first factor.
\end{proof}

It is useful to write the formula for the last factor in equation \eqref{lastZfactor} explicitly:
\begin{equation}
Z_{\lambda;\bn^{(j,p)}}^{(k-j,p)}(\bu^{(j,p)}) =
\sum_{\bm^{(j,p)}} \prod_{\al}\left[ \frac{1}{\cQ_{\al,\tau_\al+1}}\right]^{\tilde q_{\al,\tau_\al}}
\cQ_{\al,\tau_\al}^{q_{\al,\tau_\al+1}} 
\cQ_{\gamma,j}^{(\sum_{l=1}^{t_\gamma'-p-1} \Delta_{\gamma',\tau_{\gamma'}+l})}F_{j,p}(\bu)
\end{equation}
where $\bm^{(j,p)}=(m_{\al,i}: i>\tau_\al)$. Note that $\tilde q_{\al,\tau_\al}=q_{\al,\tau_\al}$ except in the case where $\al=\gamma$. Here,
\begin{equation}\label{F}
F_{j,p}(\bu) = \prod_\al \prod_{i>\tau_\al} u_{\al,i}^{q_{\al,i}} a_i^{\Delta_{\al,i}} {q_{\al,i}+m_{\al,i}\choose m_{\al,i}}
\end{equation}

\subsection{Identity of power series}

First, we need a Lemma similar to Lemma \ref{seriesform} for non-simply
laced algebras.
\begin{lemma} Let $A = \C[u_1^{-1},...,u_r^{-1}]$.
If $Q_{\al,2}^{-1}\in A[[u_\al]]$ for each $\al\in I_r$, then so is
$Q_{\al,j}^{-1}$ for all $j$.
\end{lemma}
\begin{proof}
The proof is identical to Lemma \ref{seriesform} with the slight
modification that we need to consider the $Q$-system in the more general
form
$$
\frac{1}{Q_{\al,j+1}} = \frac{Q_{\al,j-1}Q_{\al,j}^{-2}}{1-Q_{\al,j}^{-2}
\prod_{\beta\neq \al} \cT_j^{(\al,\beta)}},
$$ and we note that for all $j,\al,\beta$, $\cT_j^{(\al,\beta)}\in A$ as
it is a polynomial in $Q$'s. The proof then proceeds in exactly the
same way as Lemma \ref{seriesform}.
\end{proof}

\subsubsection{Evaluation maps}
\begin{defn}\label{generaleval}
Define the evaluation map $\varphi_{j,p}$ as follows, extended by
linearity. For $0\leq p<\max(t_\al)$, the map acts as
$$
\varphi_{j,p}(a_i) = 1 \hbox{ if $i\leq \tau_{\gamma'}$},
$$
and leaves other $a_i$'s unchanged. For $p=0$, it acts as
$$
\varphi_{j,0}(u_{\al,i}) = 1 \hbox{ if $i< \tau_{\al}$},
$$
and leaves the other $u$-variables unchanged. Finally, if $p>0$,
$$
\varphi_{j,p}(u_{\al,i}) = 1 \hbox{ if $\al$ is long and $i\leq j$ or
  if $\al$ is short and $i<\tau_\al$,}
$$
leaving the other $u$-variables unchanged.
\end{defn}

\begin{lemma}
Let $0\leq p\leq t_\al-1$. Then
$$
\cQ_{\al,t_\al j + 1+p} \hbox{ is independent of $u_{\beta,i}$ for
  $i\geq t_\beta j+\delta_{p>0}$ and $\beta\neq \al$ }.
$$
and
$$
\cQ_{\al,j} \hbox{ is independent of $u_{\al,i}$ with $i\geq j$.}
$$
\end{lemma}

\begin{proof}
First, note that for each $j$, $\cQ_{\al,j+1}$ is a function of
$u_{\al,j}$,  $\cQ_{\al,i}$ with $i\leq j$ and
$\cT_j^{(\al,\beta)}$ with $\beta\sim \al$.

Suppose the statement of the lemma is true for all $l<j$ and all $p$,
and for $j$ with $p=0$. We proceed in three steps.

\begin{enumerate}
\item If $t_\alpha = 3$ (i.e. $\g=G_2$) then we have that (for $p=1,2,3$)
$\cQ_{2,3j+1+p}$ depends on $\cT_{3j+p}^{(2,1)}$ which depends on
$\cQ_{1,i}$ with $i\leq j+1$. We can use the induction hypothesis (for
$t_\al=1$) to deduce that this is independent of $u_{\beta,i}$ with $i\geq
t_\beta j + 1$ for $\beta=1,2$. Moreover for each $p$, we have an
explicit dependence on $u_{2,3j+p}$. By induction on $p$, we conclude that
$\cQ_{2,3j+1+p}$ is independent of $u_{1,i}$ with $i>j$ and $u_{2,i}$ with $i\geq 3j+1+p$.

For $t_\alpha=3$, we conclude that the statement of the Lemma holds
for $j$ with $p=1,2$ and for $j+1$ with $p=0$.

\item If $t_\al=2$ consider first $\cQ_{\al,2j+2}$ which depends on
$\cT_{2j+1}^{(\al,\beta)}$, which depends on $\cQ_{\beta,i}$ with $i\leq
t_\beta j + 1$. By the induction hypothesis this is independent of
$u_{\delta,i}$ with $i\geq t_\delta j + 1$ for any root
$\delta$. Thus, $\cQ_{\alpha, 2j+2}$ is independent of
$u_{\beta,2j+1}$ if $\beta\neq \al$ and of $u_{\al,2j+2}$.
We conclude that the statement of the Lemma holds for $j$ with $p=1$. 

We then consider $\cQ_{\al,2j+3}$ which depends on
$\cT_{2j+2}^{(\al,\beta)}$. This depends on
$\cQ_{\beta,t_\beta(j+1)}$, with $t_\beta = 1$ or $t_\beta = 2$. In
either case, either with the induction hypothesis or the previous
paragraph, this is independent of $u_{\delta, i}$ with $i\geq t_\delta(j+1)$.
We conclude that $\cQ_{\al,2(j+1)+1}$ is independent of $u_{\beta,i}$
with $i\geq t_\beta (j+1)$ and of $u_{\al,2j+3}$. This is the
statement of the Lemma with $p=0$ and $j+1$.

\item If $t_\al =1$ then we consider $\cQ_{\al,j+2}$. This depends on
$\cT_{j+1}^{(\al,\beta)}$ which depends on $\cQ_{\beta,
t_\beta(j+1)}$. If $t_\beta=1$ we can use the induction hypothesis. If
$t_\beta>1$ we use steps (1) or (2). In either case we can say that it
is independent of $u_{\delta,i}$ with $i\geq t_\delta (j+1)$. Thus, we
have the statement of the Lemma for $j+1$.
\end{enumerate}
The Lemma follows by induction.
\end{proof}

We also have by a simple induction
\begin{lemma}
$$
 \cQ_{\al,j} \hbox{ is independent of $a_i$ with $i\geq t_{\gamma'}j/t_\al$.}
$$
\end{lemma}
We have the easy corollary:
\begin{cor}\label{evalcor}
\begin{eqnarray*}
\varphi_{j,0} (\cQ_{\al,i}) &=& Q_{\al,i}\ (i\leq t_\al j) \\
\varphi_{j,0} (\cQ_{\al, t_\al j + 1})  &=& \frac{Q_{\al,t_\al
j+1}}{u_{\al,t_\al j}},\\
\varphi_{j,p} (\cQ_{\al,i}) &= & Q_{\al,i} \ (p>0,\ i\leq \tau_\al +
\delta_{t_\al,1}), \\
\varphi_{j,p}( \cQ_{\al,\tau_\al+1} )&= &  \frac{Q_{\al,\tau_\al
+1}}{u_{\al,\tau_\al}},\ \hbox{($\al$ short)}.
\end{eqnarray*}
\end{cor}
\subsubsection{Series expansions}
\begin{defn}
Let $K=\{k_1,...,k_r\}\in \N^r$. The restricted summation
$Z_{\lambda;\bn}^{(k)}(\bu)^{[K]}$ is the generating function
\eqref{zgen} restricted to $\bm$ such that $q_{\al,j}\geq 0$ for all
$k_\al \leq j <k$ and for each $\al\in I_r$. 
\end{defn}

\begin{lemma}\label{first}
For fixed $j$ and $\al$, define $K=\{k_1, ..., k_r \}$ with
$k_\al=t_\al j+1$ and $k_\beta\geq t_\beta j$ ($\beta\in I_r$). Then
$$
\cP_{u_\al}
\varphi_{j,0}(Z_{\lambda;\bn}^{(k)}(\bu)^{[K]})
=
\cP_{u_\al}
\varphi_{j,0}(Z_{\lambda;\bn}^{(k)}(\bu)^{[K-\epsilon_\al]})
$$
\end{lemma}
\begin{proof}
Consider the factorization
$$
Z_{\lambda;\bn}^{(k)}(\bu) = Z_{0;\bn}^{(j)}(\bu) \
Z_{\lambda;\bn^{(j)}}^{(k-j)}(\bu^{(j)}) .
$$
Let $K=(k_1,...,k_r)$ with $k_\al \geq t_\al j$, and consider the restricted generating function. The factorization formula still holds, with the restriction affecting only the second factor, because the first factor does not depend on the variables $u_{\beta,i}$ with $i\geq \tau_\beta$ for each $\beta$. Therefore,
$$
Z_{\lambda;\bn}^{(k)}(\bu)^{[K]} = Z_{0;\bn}^{(j)}(\bu) \
Z_{\lambda;\bn^{(j)}}^{(k-j)}(\bu^{(j)})^{[K]}  .
$$
 In the first
factor, we have a product of $\cQ_{\al,i}$ with $i\leq t_\al j +
1$. Applying the evaluation map $\varphi_{j,0}$ and using Corollary \ref{evalcor},
the first factor becomes (after a cancellation of the factors
$u_{\al,t_\al j}$)
$$
\varphi_{j,0} (Z_{0;\bn}^{(j)}(\bu)) = \prod_{\al} \frac{Q_{\al,1}
  Q_{\al,t_\al j}}{Q_{\al,t_\al j + 1}} \prod_{i=1}^{t_\al j}
Q_{\al,i}^{n_{\al,i}}. 
$$
The second factor can be written out explicitly as a summation over
the variables $\bm^{(j)}:=\{m_{\al,i}| \ \al\in I_r, \ i>t_\al j\}$:
$$
Z_{\lambda;\bn^{(j)}}^{(k-j)} (\bu^{(j)})^{[K]}  =
\sum_{\bm^{(j)}\atop q_{\al,i}\geq 0 \ (i\geq k_\al)} \prod_\al \left[\frac{1}{\cQ_{\al, t_\al j +
    1}(\bu)}\right]^{q_{\al, t_\al j}} \cQ_{\al,t_\al j}(\bu)^{q_{\al,t_\al j+1}}
F_{j,0}(\bu).
$$
Here, $F_{j,0}(\bu)$ is defined in equation \eqref{F}.
The evaluation map $\varphi_{j,0}$ has the following effect only:
$$
\varphi_{j,0}(Z_{\lambda;\bn^{(j)}}^{(k-j)} (\bu^{(j)})^{[K]} ) =
\sum_{\bm^{(j)}\atop q_{\al,i}\geq 0 \ (i\geq k_\al)} \prod_\al \left[\frac{u_{\al,t_\al j}}{Q_{\al, t_\al j +
    1}}\right]^{q_{\al, t_\al j}} Q_{\al,t_\al j}^{q_{\al,t_\al j+1}}
F_{j,0}(\bu).
$$
For a fixed $\al$, we are interested in non-negative powers in $u_\al$ in the Laurent expansion of this function in $u_\al$. The dependence on $u_\beta$ for all $\beta$ is contained in the factors involving $Q_{\beta,i}$. Moreover the Laurent expansion has positive powers of $u_\al$ which come only from the factors with $\beta=\al$.

Assume now that $k_\al>t_\al j$ for some $\al$. Consider a term in the summation over $\bm^{(j)}$ such that
$q_{\al, t_\al j + 1}\geq 0$ is fixed and $q_{\al,t_\al j}<0$ for some $\al$. Considering the dependence on the factors $Q_{\al,i}$ for various $i$ we see that there are no terms left in the denominator as in this case the denominator $Q_{\al,\tau_\al+1}$ in the prefactor cancels and all other terms have non-negative powers. Therefore we have a polynomial in $Q_{\al,i}$ and hence in $u_\al^{-1}$. Moreover, this polynomial has no constant term in $u_\al^{-1}$ because there is an overall factor of $Q_{\al,1}=u_{\al}^{-1}$.
Thus, the power series
expansion in $u_\al$ has no contribution from this term.
The power series in $u_\al$ of the term with $q_{\al,\tau_\al+1}\geq 0$
 has nontrivial contributions only from terms with $q_{\al,\tau_\al}\geq 0$.
\end{proof}

\begin{lemma}\label{second}
For $0<p<\max(t_\al)$, and for $\al$ a fixed short root, let
$K=\{k_1,...,k_r\}$ with $k_\beta \geq \tau_\beta$ for all $\beta$ and $k_\al = \tau_\al+1$. Then
$$
\cP_{u_\al}\varphi_{j,p}(Z_{\lambda;\bn}^{(k)}(\bu)^{[K]}) =
\cP_{u_\al}\varphi_{j,p}(Z_{\lambda;\bn}^{(k)}(\bu)^{[K-\epsilon_\al]}).
$$
\end{lemma}
\begin{proof}
Consider the factorization \eqref{partialfactorization} for
some fixed $j$ and $p$. Fix $K=(k_1,...,k_r)$ with $k_\beta\geq \tau_\beta$. The factorization \eqref{partialfactorization} still holds for the restricted generating function $Z_{\lambda;\bn}^{(k)}(\bu)^{[K]}$ because the restriction only affects the second factor -- the first factor does not depend on the variables $u_{\al,i}$ with $i\geq \tau_\al$.
Apply the evaluation map
$\varphi_{j,p}$ to the factorized formula. Using corollary \ref{evalcor}, we have 
$$
\varphi_{j,p}(Z_{\lambda;\bn}^{(k)}(\bu)^{[K]}) =
\prod_{\al} \frac{Q_{\al,1}Q_{\al,\tau_\al}} {Q_{\al,\tau_\al+1}}
\prod_{i=1}^{\tau_\al} {Q_{\al,i}^{n_{\al,i}}}\times \varphi_{j,p}(
Z_{\lambda;\bn^{(j,p)}}^{(k-j,p)}(\bu^{(j,p)})^{[K]})
$$
The second factor is
$$
\sum_{m_{\al,i}\ (i>\tau_\al)\atop q_{\al,i}\geq k_\al}
\left[\prod_{\al\in \Pi^>} Q_{\al,j+1}^{-\tilde{q}_{\al,j}}\right]
\left[\prod_{\al\in \Pi^<}
\frac{u_{\al,\tau_\al}}{Q_{\al,\tau_\al+1}}\right]^{{q}_{\al,\tau_\al}}
\prod_{\al\in I_r} (u_{\al,\tau_\al+1}
Q_{\al,\tau_\al})^{q_{\al,\tau_\al+1}} Q_{\gamma,j}^\mu F_{j,p}(\bu).
$$
Here, $\mu=\sum_{i=1}^{t_{\gamma'}-p-1}\Delta_{\gamma',\tau_{\gamma'}+i}\geq 0$.
Note that if $i> \tau_\al$, then
\begin{equation}\label{inequality}
\tilde q_{\beta,i}:=q_{\beta,i}|_{m_{\al,\ell}=0\ (\ell\leq
  \tau_\al)}=q_{\beta,i}-\delta_{\beta,\gamma}\delta_{i,j} 
\sum_{\ell=1}^p \ell \ m_{\gamma',\tau_{\gamma'}-p+\ell}.
\end{equation}

Consider a term in the summation over $\bm^{(j,p)}$ with fixed 
$q_{\al,\tau_\al+1}\geq 0$ and $q_{\al,\tau_\al}<0$ for some short root $\al$.
The positive powers in $u_\al$ in the generating function can come only from the expansion of the functions $Q_{\al,i}$ for various $i$ appearing in the factorized formula above. Examination of the dependence on such factors in this case shows that there are no negative powers of $Q_{\al,i}$, and hence such a term is a polynomial in $u_\al^{-1}$ with no constant term due to the presence of the factor $Q_{\al,1}=u_\al^{-1}$ in the prefactor.
Therefore, if $q_{\al, \tau_\al+1}\geq 0$, the power series expansion in $u_\al$ has contributions only from terms with $q_{\al,\tau_\al}\geq 0$.
\end{proof}

Finally we note that if we replace the evaluation maps $\varphi_{j,p}$
which appear in Lemmas \ref{first} and \ref{second} with
$\varphi_{k+1,0}$ (that is, evaluation at $u_{\beta,i}=1=a_i$ for all $i,\beta$) the Lemmas still hold.


From this follows our main theorem, which implies the $M=N$ identity
for non-simply laced algebras.

\begin{thm}\label{main}
Let $K=(1,...,1)$. Then
\begin{equation}\label{finalequation}
\cP_{u_1,...,u_r} \varphi_{k+1,0}(Z_{\lambda;\bn}^{(k)}(\bu)) =
\cP_{u_1,...,u_r} \varphi_{k+1,0}(Z_{\lambda;\bn}^{(k)}(\bu)^{[K]}) .
\end{equation}
\end{thm}
\begin{proof}
Lemmas \ref{first} and \ref{second} can be applied successively to several roots
$\al_1,...,\al_t$. Thus we can write, when $t=r$,
$$
\cP_{u_1,...,u_r}\varphi_{j+1,0}(Z_{\lambda;\bn}^{(k)}(\bu)^{[K]})=
\cP_{u_1,...,u_r}\varphi_{j+1,0}(Z_{\lambda;\bn}^{(k)}(\bu)^{[K']}),
$$
where $K$ is the set with $k_\al=t_\al (j+1)$ and $K'$ is the set with
$k_\al=t_\al j$ for all $\al\in I_r$. The Theorem follows by induction
on $j$.
\end{proof}

\begin{cor} 
The identity \eqref{mainidentity} is true for any simple Lie algebra $\g$.
\end{cor}
\begin{proof} As for the simply-laced Lie algebras in the previous section, 
this is just the constant term of Equation \eqref{finalequation}.
\end{proof}

\section{Conclusion}
In this paper we have proved a combinatorial identity, Conjecture
\ref{HKOTYconj}, which implies the fermionic form $M_{\lambda;\bn}$
for the multiplicities of the irreducible $\g$-modules in the tensor
product of KR-modules for untwisted Yangians or quantum affine algebras. This was done in each case by constructing an appropriate generating function satisfying a factorization property, which allows us to prove that the restricted $M$-sum is equal to the unrestricted $N$-sum. This method appears to be quite general.

There are certain generalizations of these formulas for twisted Yangians
\cite{twistedHKOTY}, and the Kirillov-Reshetikhin characters
in these cases were shown to satisfy a $Q$-system \cite{He07}. We have
not addressed these systems in this paper, although it is clear that
the same methods should apply to these cases. More generally, it would
be interesting to understand what is the most general form of vacancy
numbers, allowing for the exact cancellations resulting in an
$M=N$-type identity.

The main structures introduced in this paper, the factorizing
generating functions and the deformed $\cQ$-systems, clearly require
further study. For example, it would be interesting to understand them
from a representation-theoretical point of view.

The most important property of $Q$-systems we used in this paper is
their polynomiality, Theorem \ref{polynomiality}. This follows from
representation theory. However, $Q$-systems can also be expressed in
the context of cluster algebras \cite{Ke07}, where a similar property
known as the Laurent phenomenon is satisfied. Polynomiality is a very
special subcase of this phenomenon, which awaits further study.

\vskip.2in
\noindent{\bf Acknowledgements:} 
We thank V. Chari, N. Reshetikhin, D. Hernandez for their
valuable input. We also thank the referee for his careful reading of
the manuscript and helpful remarks. RK thanks CTQM at the University
of Aarhus and of CEA-Saclay IPhT for their hospitality.  RK is
supported by NSF grant DMS-05-00759. PDF acknowledges the support of
European Marie Curie Research Training Networks ENIGMA
MRT-CT-2004-5652, ENRAGE MRTN-CT-2004-005616, ESF program MISGAM, ACI
GEOCOMP and of ANR program GIMP ANR-05-BLAN-0029-01.

\begin{appendix}
\section{Proof of Lemmas \ref{partialfactorization}, \ref{lemmageninit}, and \ref{recuZ}}\label{Zone}
In this appendix, we prove the summation Lemmas \ref{partialfactorization}, \ref{lemmageninit}, and \ref{recuZ}.
To do so, we start from $Z^{(k)}_{\lambda;\bn}(\bu)$ as defined through \eqref{zgen},
and explicitly sum over some $m_{\al,i}$, with $i\leq p \leq t_\al$, in a specific order.

\subsection{Preliminaries}

We recall the definition \eqref{zgen} of the generating fuction:
$$
Z_{\lambda;\bn}^{(k)}(\bu) = \sum_{\bm}
\prod_{\al=1}^r u_\al^{q_\al} \prod_{i=1}^{t_\al k}
u_{\al,i}^{q_{\al,i}} a_i^{\Delta_{\al,i}}
{m_{\al,i}+q_{\al,i}\choose m_{\al,i}} ,
$$
where we list the modified vacancy numbers for each $\g$ explicitly, 
from the definition \eqref{nonsimpleq} (note that $q_{\al,0}=q_\al$):
\begin{eqnarray*}
B_r: &&q_{\al,i}=l_\al+\sum_{j=i+1}^k
   (j-i)(2m_{\al,j}-n_{\al,j}-m_{\al-1,j}-m_{\al+1,j}),\quad (\al<r-1)
\\ && q_{r-1,i}=
   l_{r-1}+\sum_{j=i+1}^k (j-i)(2m_{r-1,j}-n_{r-1,j}-m_{r-2,j})
   -\sum_{j=2i+1}^{2k} (j-2i)m_{r,j},\\ &&
   q_{r,i}=l_r+\sum_{j=i+1}^{2k}
   (j-i)(2m_{r,j}-n_{r,j})-\sum_{i<2j\leq 2k} (2j-i)m_{r-1,j}.
\end{eqnarray*}
\begin{eqnarray*}
C_r: && q_{\al,i}=l_\al+\sum_{j=i+1}^{2k}
  (j-i)(2m_{\al,j}-n_{\al,j}-m_{\al-1,j}-m_{\al+1,j}),\quad (\al<r-1)\\ &&
  q_{r-1,i}=l_{r-1}+\sum_{j=i+1}^{2k} (j-i)(2m_{r-1,j}-n_{r-1,j})
  -\sum_{j=i+1}^{2k} (j-i)m_{r-2,j}\\ 
&& \hskip1in-\sum_{i<2j\leq 2k} (2j-i)m_{r,j},\\
  && q_{r,i}=l_r+\sum_{j=i+1}^k
  (j-i)(2m_{r,j}-n_{r,j})-\sum_{j=2i+1}^{2k} (j-2i)m_{r-1,j}.
\end{eqnarray*}
\begin{eqnarray*}
F_4: && q_{1,i}=l_1+\sum_{j=i+1}^k (j-i)(2m_{1,j}-n_{1,j}-m_{2,j})\\
  && q_{2,i}=l_2+\sum_{j=i+1}^k (j-i)(2m_{2,j}-n_{2,j}-m_{1,j})
-\sum_{j=2i+1}^{2k} (j-2i)m_{3,j}\\
  && q_{3,i}=l_3+\sum_{j=i+1}^{2k} (j-i)(2m_{3,j}-n_{3,j}-m_{4,j})
-\sum_{i<2j\leq 2k} (2j-i)m_{2,j} \\
  && q_{4,i}=l_4+\sum_{j=i+1}^{2k} (j-i)(2m_{4,j}-n_{4,j}-m_{3,j})
\end{eqnarray*}
\begin{eqnarray*}
G_2: &&  q_{1,i}=l_1+
\sum_{k\geq j>i} (j-i)(2m_{1,j} -n_{1,j}) -\sum_{3k\geq j>3i} (j-3i)m_{2,j}\\
  && q_{2,i}=l_2+
\sum_{3k\geq j>i} (j-i)(2m_{2,j} -n_{2,j} )-\sum_{3k\geq 3j>i} (3j-i)m_{1,j}.
\end{eqnarray*}
Here, we use the convention that $m_{0,i}:=0$.
As before, we use the notation $\gamma'$ and $\gamma$ for the short 
and long roots connected to each other in the Dynkin diagram, namely
$\gamma'=r,r-1,3,2$ and $\gamma=r-1,r,2,1$ for $B_r,C_r,F_4,G_2$ 
respectively.

We will sum \eqref{zgen} over the $m_{\al,i}$,
$i=1,...,t_\al$ and $\al=1,...,r$ explicitly. The summation over these integers must be done in a certain  order, because in each case, the $\{ q_{\gamma',i} |
i=1,...,t_{\gamma'}-1\}$ depend on $m_{\gamma,1}$.  We must
therefore first sum over $m_{\gamma',i}$, $i=1,...,t_{\gamma'}-1$ before summing over the other variables.
 The intermediate summations are
rational fractions of the $u$'s and $a$'s, which can be expressed
in terms of the functions $\cQ_{\al,i}$.  

We first present the partial summations leading to Lemma \ref{partialfactorization},
for which only the $m_{\beta,i}$, $\beta \in \Pi^<$, $i=1,2,...,p<t_{\gamma'}$ are summed
over, starting from $\beta=\gamma'$,
and then the ``complete" sums corresponding to $p=t_{\gamma'}$
for short root $m$'s and also summed over the long root $m$'s, leading to Lemmas \ref{lemmageninit} and \ref{recuZ}.

A crucial ingredient used repeatedly in the following is the fact that the $q$'s
satisfy relations expressing $q_{\al,i-1}$ in terms
solely of $m_{\al,i}, n_{\al,i}$ and the combination $2q_{\al,i}-q_{\al,i+1}$
with possible slight modifications, involving only finitely many $m$'s. 
These relations are (we use the notation $q_{\al,0}:= q_\al$
for convenience):
\begin{eqnarray}\label{qs}
&&\bullet \ {\rm short \ root}\ \gamma': \nonumber\\
&&q_{\gamma',j-1}=-n_{\gamma',j}+ \sum_{\beta\in \Pi^<} C_{\gamma',\beta} m_{\beta,j}
+2 q_{\gamma',j}-q_{\gamma',j+1}-\delta_{j,t_{\gamma'}}m_{\gamma,1}, \quad 1\leq j\leq t_{\gamma'} 
\label{qgamp}\\
&&\bullet \ {\rm short \ root}\ \beta\neq \gamma': \nonumber\\
&&q_{\beta,j-1}=-n_{\beta,j}+ \sum_{\beta'\in \Pi^<} C_{\beta,\beta'} m_{\beta',j}
+2 q_{\beta,j}-q_{\beta,j+1} \label{qbeta}\\
&&\bullet \ {\rm long \ root}\ \gamma: \nonumber\\
&&q_{\gamma}= -n_{\gamma,1}+\sum_{\beta\in \Pi^>} C_{\gamma,\beta} m_{\beta,j}
+2q_{\gamma,1}-q_{\gamma,2}- \sum_{j=1}^{t_\gamma'}( j m_{\gamma',j}
+\Delta_{\gamma',t_{\gamma'}+j})\label{qgam}\\
&&\bullet \ {\rm long \ root}\ \alpha\neq \gamma: \nonumber\\
&&q_{\al}= -n_{\al,1}+\sum_{\beta\in \Pi^>} C_{\al,\beta} m_{\beta,j}
+2q_{\al,1}-q_{\al,2} \label{qal}
\end{eqnarray}

\subsection{Partial summations over short roots: 
proof of Lemma \ref{partialfactorization}}\label{apptwo}

In each case, we must first sum over $\mu=m_{\gamma',1}$. Collecting
all relevant factors in the summand of \eqref{zgen} and terms which
depend on $u_{\gamma'}$ and $u_{\gamma',1}$, using \eqref{qgamp} for
$j=1$, we have
$$ \sum_{\mu\geq 0} u_{\gamma'}^{-n_{\gamma',1}+2
q_{\gamma',1}-q_{\gamma',2}} \left( a_1^{t_{\gamma'}-1} \prod_\beta
u_\beta^{C_{\beta,\gamma'}} \right)^{\mu}
u_{\gamma',1}^{q_{\gamma',1}} {\mu+q_{\gamma',1}\choose \mu}
={\cQ_{\gamma',1}^{n_{\gamma',1}}\over u_{\gamma',1}}
{\cQ_{\gamma',1}^2 \over \cQ_{\gamma',2}}
\frac{\cQ_{\gamma',1}^{q_{\gamma',2}}}{\cQ_{\gamma',2}^{q_{\gamma',1}}
} ,
$$
where we have identified $\cQ_{\gamma',1}=u_{\gamma'}^{-1}$ and
$$\cQ_{\gamma',2}={(1-a_1^{t_{\gamma'}-1}u_{\gamma}^{-1}\prod_{\beta\in \Pi^<} u_\beta^{C_{\beta,\gamma'}})u_{\gamma'}^{-2}\over u_{\gamma',1}}.$$

Assume $t_{\gamma'}=2$. Then we may now sum over $m_{\beta,1}$
for the other short roots $\beta\neq \gamma'$, where
 $t_\beta=2$. Using the relation \eqref{qbeta} for $j=1$, 
we have for each $\beta\in \Pi^>$, $\beta\neq \gamma'$, a factor of the form
$$
\sum_{\mu\geq 0}
u_{\beta}^{-n_{\beta,1}+2 q_{\beta,1}-q_{\beta,2}} \left(\prod_{\beta'}
u_{\beta'}^{C_{\beta',\beta}} \right)^{\mu} u_{\beta,1}^{q_{\beta,1}}
{\mu+q_{\beta,1}\choose \mu} ={\cQ_{\beta,1}^{n_{\beta,1}}\over u_{\beta,1}} 
{\cQ_{\beta,1}^2 \over \cQ_{\beta,2}}  
\frac{\cQ_{\beta,1}^{q_{\beta,2}}}{\cQ_{\beta,2}^{q_{\beta,1}} } ,
$$
where we have identified $\cQ_{\beta,1}=u_{\beta}^{-1}$ and
$$\cQ_{\beta,2}={(1-\prod_{\beta'} u_{\beta'}^{C_{\beta',\beta}})u_{\beta}^{-2}\over u_{\beta,1}}.$$
Gathering all the above contributions and restricting 
$q_\gamma$ to ${\tilde q}_\gamma\equiv q_\gamma\vert_{m_{\gamma',1}=0}$ yields
Lemma \ref{partialfactorization}  for $p=1$, the substitutions \eqref{subszerop} being induced by the factors $\cQ_{\beta,1}^{q_{\beta,2}}/\cQ_{\beta,2}^{q_{\beta,1}}$, $\beta\in \Pi^<$.

If $t_{\gamma'}=3$ ($G_2$ case, $\gamma'=2$, $\gamma=1$) then if $p=1$, we have proved Lemma \ref{partialfactorization} above. If $p=2$,
we must now sum over $\mu=m_{2,2}$. We use \eqref{qgamp} for $j=2$
to rewrite $q_{\gamma',1}=2\mu-n_{\gamma',2}+2q_{\gamma',2}-q_{\gamma',3}$ in
the summation:
\begin{equation}\label{gtwopart}
{\cQ_{\gamma',1}^{n_{\gamma',1}}\over u_{\gamma',1}} 
{\cQ_{\gamma',1}^2 \over \cQ_{\gamma',2}}  \sum_{\mu\geq 0}
\frac{\cQ_{\gamma',1}^{q_{\gamma',2}}}{\cQ_{\gamma',2}^{q_{\gamma',1}}} a_2^{\mu} u_\gamma^{-2\mu}
u_{\gamma',2}^{q_{\gamma',2}}
{\mu+q_{\gamma',2}\choose \mu} 
={\cQ_{\gamma',1}^{n_{\gamma',1}}\over u_{\gamma',1}} 
{\cQ_{\gamma',2}^{n_{\gamma',2}}\over u_{\gamma',2}} \cQ_{\gamma',1}{\cQ_{\gamma',2}\over \cQ_{\gamma',3}}  
\frac{\cQ_{\gamma',2}^{q_{\gamma',3}}}{\cQ_{\gamma',3}^{q_{\gamma',2}} } 
\end{equation}
where we have identified 
$$\cQ_{\gamma',3}={\cQ_{\gamma',2}^2-a_2 \cQ_{\gamma,1}^2\over u_{\gamma',2}\cQ_{\gamma',1}}.$$
Replacing in the summation formula for $Z^{(k)}_{\lambda;\bn}(\bu)$ the quantity $q_\gamma$ by 
${\tilde q}_{\gamma}\equiv q_{\gamma}\vert_{m_{\gamma',1}=m_{\gamma',2}=0}$ 
yields the Lemma \ref{partialfactorization} for $p=2$ ($G_2$ case).

\subsection{Complete summations: proof of Lemma \ref{lemmageninit}}\label{appthree}

Starting from Lemma \ref{partialfactorization} with $p=t_{\gamma'}-1$, we have two cases to consider.
In the case $t_{\gamma'}=2$ ($B_r,C_r,F_4$), let us first sum over 
$\mu=m_{\gamma',2}$. From \eqref{qgamp} 
for $j=2=t_{\gamma'}$, we get $q_{\gamma',1}=2\mu-n_{\gamma',2}+2q_{\gamma',2}-q_{\gamma',3}
-m_{\gamma,1}$, hence the summation:
$$
\sum_{\mu\geq 0} {\cQ_{\gamma',1}^{n_{\gamma',1}}\over u_{\gamma',1}} 
{\cQ_{\gamma',1}^2 \over \cQ_{\gamma',2}}  
\frac{\cQ_{\gamma',1}^{q_{\gamma',2}}}{\cQ_{\gamma',2}^{q_{\gamma',1}}}  u_\gamma^{-2\mu}
u_{\gamma',2}^{q_{\gamma',2}}
{\mu+q_{\gamma',2}\choose \mu} ={\cQ_{\gamma',1}^{n_{\gamma',1}}\over u_{\gamma',1}} 
{\cQ_{\gamma',2}^{n_{\gamma',2}}\over u_{\gamma',2}} \cQ_{\gamma',1}{\cQ_{\gamma',2}\over \cQ_{\gamma',3}}  
\frac{\cQ_{\gamma',2}^{q_{\gamma',3}+m_{\gamma,1}}}{\cQ_{\gamma',3}^{q_{\gamma',2}} } 
$$
where  
$$\cQ_{\gamma',3}={\cQ_{\gamma',2}^2-\cQ_{\gamma,1}^2\over u_{\gamma',2}\cQ_{\gamma',1}}.$$
We may now sum on $\mu=m_{\beta,2}$ for the other short roots $\beta\neq \gamma'$.
We use \eqref{qbeta} for $j=2$, and compute:
$$
\sum_{\mu\geq 0} {\cQ_{\beta,1}^{n_{\beta,1}}\over u_{\beta,1}} 
{\cQ_{\beta,1}^2 \over \cQ_{\beta,2}}  
\frac{\cQ_{\beta,1}^{q_{\beta,2}}}{\cQ_{\beta,2}^{q_{\beta,1}}} 
u_{\beta,2}^{q_{\beta,2}}
{\mu+q_{\beta,2}\choose \mu} ={\cQ_{\beta,1}^{n_{\beta,1}}\over u_{\beta,1}} 
{\cQ_{\beta2}^{n_{\beta,2}}\over u_{\beta,2}} \cQ_{\beta,1}{\cQ_{\beta,2}\over \cQ_{\beta,3}}  
\frac{\cQ_{\beta,2}^{q_{\beta,3}}}{\cQ_{\beta,3}^{q_{\beta,2}} } 
$$
where 
$$\cQ_{\beta,3}={(1-\prod_{\beta'} Q_{\beta',2}^{-C_{\beta',\beta}})
\cQ_{\beta,2}^2\over u_{\gamma',2}\cQ_{\gamma',1}}.$$
We may now finally sum over $m_{\al,1}$ for all $\al \in \Pi^>$. We first do the summation
over $\mu=m_{\gamma,1}$, using \eqref{qgam}:
$$
\sum_{\mu\geq 0} u_{\gamma}^{-n_{\gamma,1}-\Delta_{\gamma',3}} \left( \cQ_{\gamma',2} \prod_{\al\in \Pi^>}
u_\al^{C_{\al,\gamma}} \right)^\mu u_{\gamma,1}^{q_{\gamma,1}}
{\mu+q_{\gamma,1}\choose \mu} =
{\cQ_{\gamma,1}^{n_{\gamma,1}}\over u_{\gamma,1}}  {\cQ_{\gamma,1}^2\over \cQ_{\gamma,2}}  
\cQ_{\gamma,1}^{\Delta_{\gamma',3}}
\frac{\cQ_{\gamma,1}^{q_{\gamma,2}}}{\cQ_{\gamma,2}^{q_{\gamma,1}} } 
$$
where we have identified $\cQ_{\gamma,1}=u_\gamma^{-1}$ and
$\cQ_{\gamma,2}=(1-\cQ_{\gamma',2}\prod_{\al\in \Pi^>}
u_\al^{C_{\al,\gamma}})u_\gamma^{-2}/u_{\gamma,1}$.
Next, we do the remaining summations over each $\mu=m_{\al,1}$, for long roots
$\al\neq \gamma$, and use \eqref{qal}:
$$
\sum_{\mu\geq 0} u_\al^{-n_{\al,1}+2q_{\al,1}-q_{\al,2}} 
\left( \prod_{\al'} u_{\al'}^{C_{\al',\al}} \right)^\mu u_{\al,1}^{q_{\al,1}}
{\mu+q_{\al,1}\choose \mu}={\cQ_{\al,1}^{n_{\al,1}}\over u_{\al,1}}  {\cQ_{\al,1}^2\over \cQ_{\al,2}} 
\frac{\cQ_{\al,1}^{q_{\al,2}}}{\cQ_{\al,2}^{q_{\al,1}} } 
$$
where  $\cQ_{\al,1}=u_\al^{-1}$ and
$$\cQ_{\al,2}={(1-\prod_{\al'} u_{\al'}^{C_{\al',\al}} )u_\al^{-2}\over u_{\al,1}}.$$

In the case $t_{\gamma'}=3$ ($G_2$), let us first sum over 
$\mu=m_{\gamma',3}$, after summing over $m_{\gamma',1},m_{\gamma',2}$
as in \eqref{gtwopart}. Using \eqref{qgamp} for $j=t_{\gamma'}=3$, namely 
$q_{\gamma',2}=2\mu -n_{\gamma',3}+2q_{\gamma',3}-q_{\gamma',4}-m_{\gamma,1}$,
this gives
$$
{\cQ_{\gamma',1}^{n_{\gamma',1}}\over u_{\gamma',1}} 
{\cQ_{\gamma',2}^{n_{\gamma',2}}\over u_{\gamma',2}} \cQ_{\gamma',1}{\cQ_{\gamma',2}\over \cQ_{\gamma',3}}  
\sum_{\mu\geq 0} u_{\gamma}^{-3\mu}
\frac{\cQ_{\gamma',2}^{q_{\gamma',3}}}{\cQ_{\gamma',3}^{q_{\gamma',2}} } 
u_{\gamma',3}^{q_{\gamma',3}} {\mu+q_{\gamma',3}\choose \mu} =
{\cQ_{\gamma',1}^{n_{\gamma',1}}\over u_{\gamma',1}} 
{\cQ_{\gamma',2}^{n_{\gamma',2}}\over u_{\gamma',2}}
{\cQ_{\gamma',3}^{n_{\gamma',3}}\over u_{\gamma',3}}
\cQ_{\gamma',1}{\cQ_{\gamma',3}\over \cQ_{\gamma',4}}  
\frac{\cQ_{\gamma',3}^{q_{\gamma',4}+m_{\gamma,1}}}{\cQ_{\gamma',4}^{q_{\gamma',3}} } 
$$
where  
$$\cQ_{\gamma',4}={\cQ_{\gamma',3}^2-u_{\gamma}^{-3}\over u_{\gamma',3}\cQ_{\gamma',2}}.$$ We next sum over $\mu=m_{\gamma,1}$, and use \eqref{qgam}. Since we have already
summed over $m_{2,i}$  $(i=1,2,3)$ these integers should be  set to zero in the definition of $q_\gamma$, namely
$\tilde q_\gamma=2\mu-n_{\gamma,1}+2q_{\gamma,1}-q_{\gamma,2}
-\Delta_{\gamma',4}-\Delta_{\gamma',5}$,
to finally write:
$$
u_\gamma^{-n_{\gamma,1}} \sum_{\mu\geq 0} \left( \cQ_{\gamma',3} u_{\gamma}^2 \right)^\mu 
u_{\gamma,1}^{q_{\gamma,1}}
{\mu+q_{\gamma,1}\choose \mu}={\cQ_{\gamma,1}^{n_{\gamma,1}}\over u_{\gamma,1}} 
\frac{\cQ_{\gamma,1}^2}{\cQ_{\gamma,2}} 
\cQ_{\gamma,1}^{\Delta_{\gamma',4}+\Delta_{\gamma',5}}
\frac{\cQ_{\gamma,1}^{q_{\gamma,2}}}{\cQ_{\gamma,2}^{q_{\gamma,1}}}
$$
where we have identified $\cQ_{\gamma,1}=u_\gamma^{-1}$ and
$\cQ_{\gamma,2}=(u_\gamma^{-2}-\cQ_{\gamma',3})/u_{\gamma,1}$.

Collecting all the above terms in both cases $t_{\gamma'}=2$ and $3$,
and renaming the summation
variables $m'_{\al,i}=m_{\al,i+t_\al}$, as well as $n_{\al,i}'=n_{\al,i+t_\al}$, yields in general
\begin{eqnarray}
&&Z_{\lambda;\bn}^{(k)}(\bu)=\left(\prod_{\al=1}^r
\frac{\cQ_{\al,1}\cQ_{\al,t_\al}} {\cQ_{\al,t_{\al}+1}}
\prod_{i=1}^{t_\al}
\frac{\cQ_{\al,i}^{n_{\al,i}}}{u_{\al,i}}\right)\times \label{bidule}
\\ &&\sum_{\{m'_{\al,i}\geq 0\ | (\al,i) \in J^{(k-1)}_\g\}}
\cQ_{\gamma,1}^{\mu'_\g}\prod_{\al=1}^r 
\frac{\cQ_{\al,t_\al}^{q'_{\al,1}}}
{\cQ_{\al,t_\al+1}^{q'_\al}}
\, \prod_{i=1}^{t_\al
(k-1)} u_{\al,i+t_\al}^{q'_{\al,i}} a_{i+t_{\al}}^{\Delta_{\al,i}'}
{m'_{\al,i}+q'_{\al,i}\choose
m'_{\al,i}}. \nonumber
\end{eqnarray}
Here, all the primed functions of $m$ are
the functions with $m'$ substituted for $m$, and similarly for $n$. More precisely, we have
$q'_{\al,i}=q_{\al,i+t_\al}$, $q'_\al=q_{\al,t_\al}$, and
$\mu'_\g=\sum_{j=1}^{t_{\gamma'}} \Delta_{\gamma',j}'=m'_{r,1},\ m_{r-1,1}',\ m_{3,1}',\ 2m_{2,1}'+m_{2,2}'$ respectively for $\g=B_n,C_n,F_4,G_2$. 

To complete the proof of Lemma
 \ref{lemmageninit}, in the case $k=1$,
we have summed over all $m$'s, and the sum on the r.h.s. of 
\eqref{bidule} is trivial, whereas
 $q'_\al=q'_{\al,1}=l_\al$. Eqn. \eqref{initgen} follows. 
\subsection{Proof of Lemma \ref{recuZ}} \label{lemmarecuZ}
The factorization
\eqref{Zrecursion} follows from \eqref{bidule}, upon identifying the first factor
with
$Z_{0;\bn}^{(1)}(\bu)$ (from Lemma \ref{lemmageninit} with all $l_\al=0$), 
and the second line with
$Z_{\lambda;\bn'}^{(k-1)}(\bu')$, with $\bu'$ as in \eqref{shiftedua}, the substitutions
in the $a$'s, $u_\al$'s and $u_{\al,1}$'s being induced respectively
by the factors $\cQ_{\gamma,1}^{\mu'_\g}$, $\cQ_{\al,t_\al+1}^{-q'_\al}$,
and $\cQ_{\al,t_\al}^{q'_{\al,1}}$.

\end{appendix}

\bibliographystyle{plain}
\def\cprime{$'$} \def\cprime{$'$} \def\cprime{$'$} \def\cprime{$'$}
  \def\cprime{$'$} \def\cprime{$'$} \def\cprime{$'$} \def\cprime{$'$}
  \def\cprime{$'$}

\end{document}